\documentclass[12pt]{amsart}
\usepackage{amsmath,enumerate,amsfonts,amssymb,color,graphicx,amsthm}
\usepackage{hyperref}

\usepackage{todonotes}
\usepackage[normalem]{ulem}

\def\RR{{\mathbb R}}

\def\SSphere{{\mathbb S}}

\def\Ric{{\rm Ric}}

\newtheorem{theorem}{Theorem}[section]
\newtheorem{proposition}[theorem]{Proposition}
\newtheorem{corollary}[theorem]{Corollary}
\newtheorem{lemma}[theorem]{Lemma}
\newtheorem{definition}[theorem]{Definition}
\newtheorem{remark}[theorem]{Remark}

\newtheorem{question}[theorem]{Question}

\numberwithin{equation}{section}

\def\vareps{\varepsilon}

\def\muGp{{\mu_\Gamma^+}}

\def\tr{\textrm{tr}}

\DeclareFontFamily{OT1}{rsfs}{}
\DeclareFontShape{OT1}{rsfs}{m}{n}{ <-7> rsfs5 <7-10> rsfs7 <10-> rsfs10}{}
\DeclareMathAlphabet{\mycal}{OT1}{rsfs}{m}{n}

\newcounter{marnote}

\newcommand\NewtonT[1]{{\stackrel{(#1)}{\mathop{T}}}}

\begin{document}
\title[Green's functions for nonlinear Yamabe problems]{Existence and uniqueness of Green's functions to nonlinear Yamabe problems}

\author{YanYan Li}
\address{YanYan Li, Department of Mathematics, Rutgers University, Hill Center, Busch Campus, 110 Frelinghuysen Road, Piscataway, NJ 08854, USA.}
\email{yyli@math.rutgers.edu.}
\author{Luc Nguyen}
\address{Luc Nguyen, Mathematical Institute and St Edmund Hall, University of Oxford, Andrew Wiles Building, Radcliffe Observatory Quarter, Woodstock Road, Oxford OX2 6GG, UK.}
\email{luc.nguyen@maths.ox.ac.uk.}

\date{}

\maketitle

\begin{abstract}
For a given finite subset $S$ of a compact Riemannian manifold $(M,g)$ whose Schouten curvature tensor belongs to a given cone, we establish a necessary and sufficient condition for the existence and uniqueness of a conformal metric on $M \setminus S$ such that each point of $S$ corresponds to an asymptotically flat end and that the Schouten tensor of the conformal metric belongs to the boundary of the given cone. As a by-product, we define a purely local notion of Ricci lower bounds for continuous metrics which are conformal to smooth metrics and prove a corresponding volume comparison theorem.

\end{abstract}

\tableofcontents

\section{Introduction} 

Let $(M^n,g)$ be a compact Riemannian manifold of dimension $n \ge 3$. It is well known that if the scalar curvature $R_g$ is positive, then the conformal Laplacian operator $-L_g = -\Delta_g + \frac{n-2}{4(n-1)}R_g$ has a unique positive Green's function $G_p \in C^\infty(M \setminus \{p\})$ with pole at a given point $p \in M$ such that 
\begin{equation}
L_g G_p = \delta_p \text{ on } M ,
	\label{Eq:LinearGreenEq}
\end{equation}
where $\delta_p$ is the Dirac measure centered at $p$. At the leading order, the singularity of $G_p$ at $p$ is the same as that of the Green's function for the Laplacian on $\RR^n$,
\[
G_p(x) 
	= \frac{1}{(n-2) |\SSphere^{n-1}|} d_g(x,p)^{-(n-2)} (1 + o(1)). 
\]
Here $d_g$ is the distance function with respect to $g$.

The purpose of the present paper is to establish the existence, the non-existence and uniqueness of (generalized) Green's functions when the conformal Laplacian in \eqref{Eq:LinearGreenEq} is replaced by other nonlinear operators arising in conformal geometry.

Let $\Ric_g$, $R_g$ and $A_g$ denote respectively the Ricci curvature, the scalar curvature and the Schouten tensor of $g$, 
\[
A_g = \frac{1}{n-2}(\Ric_g - \frac{1}{2(n-1)}R_g\,g),
\]
and let $\lambda(A_g)=(\lambda_1, \cdots, \lambda_n)$
denote the eigenvalues of
$A_g$ with respect to $g$. For a positive smooth function $u$, let $g_u = u^{\frac{4}{n-2}}g$. We have
\[
A_{g_u} = -\frac{2}{n-2} u^{-1} \nabla^2 u + \frac{2n}{(n-2)^2} u^{-2} du \otimes du - \frac{2}{(n-2)^2} u^{-2}|du|_g^2\,g + A_g.
\]

We are interested in constructing solutions to the equation
\[
\lambda(A_{g_u}) \in \partial\Gamma \text{ and } u > 0 \text{ away from a given finite number of points in $M$}
\]
where
\begin{equation}  
\Gamma\subset \Bbb R^n \mbox{ is an open convex symmetric
cone with vertex at  
the origin} 
\label{G1} 
\end{equation}
satisfying
\begin{equation}
\Gamma_n := \Big\{\lambda \in \Bbb R^n | \lambda_i > 0, 1 \leq i \leq n\Big\} \subset
\Gamma \subset \Gamma_1 := \Big\{\lambda \in \Bbb R^n | \sum^n_{i=1} \lambda_i >
0\Big\}.
\label{G2}
\end{equation}

Standard examples of such cones are the $\Gamma_k$ cones, $1 \leq k \leq n$, 
\[
\Gamma_k = \{\lambda \in \RR^n: \sigma_j(\lambda) > 0 \text{ for } 1 \leq j \leq k\},
\]
where $\sigma_k$ is the $k$-th elementary symmetric function, 
\[
\sigma_k(\lambda) = \sum_{i_1 < \cdots < i_k} \lambda_{i_1} \cdots \lambda_{i_k}.
\]

Note that, under \eqref{G1}-\eqref{G2}, there exists a function $f$ satisfying (see Proposition \ref{prop:fConstr} in Appendix \ref{Sec:SDefFc})
\begin{align}
&f\in C^\infty (\Gamma) \cap C^0 (\overline{\Gamma}) \mbox{ is homogeneous of degree one and symmetric in } \lambda_i,
\label{f1}
\\
&f>0\ \mbox{in}\ \Gamma,
\quad f = 0 \mbox{ on } \partial\Gamma,
\label{f2}
\\
&f_{\lambda_i} > 0 \ \mbox{in
} \Gamma   \  \forall 1 \leq i \leq n,
\label{f3}\\
&f\text{ is concave in $\Gamma$}.
\label{f4}
\end{align}
The partial differential relation $\lambda(A_{g_u}) \in \partial\Gamma$ can thus be re-expressed in a more familiar form
\[
f\big(\lambda(A_{g_u})\big) = 0.
\]

We adopt the following definition.

\begin{definition} \label{Def:GreenDef} Assume $m \geq 1$ and  let $p_1, \ldots, p_m$ be distinct points of $M$ and $c_1, \ldots, c_m$ be positive numbers.  A function $u \in C^{0}_{loc}(M \setminus \{p_1, \ldots, p_m\})$  is called a Green's function for $\Gamma$ with poles $p_1, \ldots, p_m$ and with strengths $c_1, \ldots, c_m$ if $u$ satisfies 
\begin{align}
&\lambda(A_{g_u}) \in \partial\Gamma \text{ and }  u > 0 \text{ in } M \setminus \{p_1, \ldots, p_m\},\label{Eq:DegEq1}\\
&\lim_{x \rightarrow p_i} d_g(x,p_i)^{n-2} u(x) = c_i,\qquad i = 1, \ldots, m.
\label{Eq:DegAs}
\end{align}
\end{definition}

In the above definition, \eqref{Eq:DegEq1} is understood in the viscosity sense -- see e.g. \cite{Li09-CPAM} for the definition. It follows that if $u$ is $C^2$, then $u$ satisfies \eqref{Eq:DegEq1} in the classical sense, and if $u \in C^{1,1}$, then $u$ satisfies \eqref{Eq:DegEq1} almost everywhere; see e.g. \cite[Lemma 2.5]{LiNgWang-C11Liouville}.

It should be clear that when $\Gamma = \Gamma_1$, the solution to \eqref{Eq:DegEq1}-\eqref{Eq:DegAs} is given uniquely as a linear combination of Green's functions for the conformal Laplacian with poles at $p_i$, namely $u = (n-2) |\SSphere^{n-1}| \sum_{i = 1}^m c_i G_{p_i}$.

It was known that when $(M,g)$ is conformal to the standard sphere and $m = 1$, there exists a unique Green's function for every given pole and strength. In the case $n = 4$ and $\Gamma = \Gamma_2$, this was proved in Chang, Gursky and Yang \cite{CGY02-AnnM} under $C^{1,1}$ regularity. For general cones in any dimension, this was proved in Li \cite{Li07-ARMA, Li09-CPAM} under $C^{0,1}$ regularity and in a joint work of the authors with Wang \cite{LiNgWang} under $C^0$ regularity. In fact, in this particular case the asymptotic condition \eqref{Eq:DegAs} is not needed -- it follows from these works that solutions to \eqref{Eq:DegEq1} satisfy \eqref{Eq:DegAs} for some positive constant $c_1$.

We note that, by \eqref{G1}-\eqref{G2}, equation \eqref{Eq:DegEq1} is degenerate elliptic. Furthermore it is not locally strictly elliptic if $\partial\Gamma_n \cap \partial\Gamma \neq \emptyset$. 

The motivation to consider Green's functions as in Definition \ref{Def:GreenDef} comes from the study of the $\sigma_k$-Yamabe problem
\begin{equation}
\sigma_k(\lambda(A_{g_u})) = 1, \quad \lambda(A_{g_u}) \in \Gamma_k  \text{ and } u > 0 \text{ in } M.
	\label{Eq:12VIII20-E1}
\end{equation}
This problem was first studied by Viaclovsky \cite{Viac00-Duke}.  An important aspect in the study of \eqref{Eq:12VIII20-E1} is to understand if the set of solutions to \eqref{Eq:12VIII20-E1} is compact, say in $C^2$, when $(M,g)$ is not conformally equivalent to the standard sphere. This compactness property of the solution set has been established when $k = 2$ and $n = 4$ \cite{CGY02-JAM}, or $(M,g)$ is locally conformally flat \cite{LiLi03}, or  $k > n/2$ \cite{GV07}, or $k = n/2$ \cite{LiNgPoorMan}. (For related works in the case $k > n/2$, see also \cite{TW09}.) The case $2 \leq k < n/2$ remains a major open problem. The role that Green's functions play in this context lies in the expectation that, under suitable conditions on $(M,g)$, appropriately rescaled blow-up solutions to \eqref{Eq:12VIII20-E1} converges along a subsequence to a Green's function for $\Gamma_k$. Whether this scenario holds for all manifolds $(M,g)$  and all $2\leq k < n/2$ remains to be understood. For this reason, we believe that understanding the existence of Green's functions as well as up-to-second-order estimates near the punctures for them (and rescaled solutions to \eqref{Eq:12VIII20-E1} which are close to some Green's function) will be extremely desirable.

As introduced in Li and Nguyen \cite{LiNgBocher}, let
\[
\mu_\Gamma^+ \text{ be the unique number such that } (-\mu_\Gamma^+, 1, \ldots, 1) \in \partial \Gamma.
\]
It is known that $\mu_\Gamma^+ \in [0,n-1]$.

For example, when $\Gamma = \Gamma_k$, $\mu_{\Gamma_k}^+ = \frac{n-k}{k}$. In particular, $\mu_{\Gamma_k}^+ > 1$ if and only if $k < \frac{n}{2}$ and $\mu_{\Gamma_k}^+ = 1$ for $k = \frac{n}{2}$. It is known that there is a distinctive difference between the cases $k > \frac{n}{2}$, $k = \frac{n}{2}$ and $k < \frac{n}{2}$, see e.g. Chang, Gursky and Yang \cite{CGY02-AnnM}, Guan, Viaclovsky and Wang \cite{G-V-W}, Viaclovsky \cite{Viac00-Duke}. Likewise, for general cones $\Gamma$, the differential inclusion $\lambda(A_g) \in \Gamma$ is sensitive to whether $\mu_\Gamma^+$ is larger, smaller or equal to $1$, see \cite{LiNgBocher}. The existence of Green's functions is also influenced by $\mu_\Gamma^+$, namely we show that they exist if and only if $\mu_\Gamma^+ > 1$, unless $(M,g)$ is conformal to the standard sphere and $m = 1$. We also prove that Green's functions, if exist, are unique. We would like to remark that the uniqueness is not straightforward, in light of the known failure of the strong maximum principle for \eqref{Eq:DegEq1}, cf. Li and Nirenberg \cite{LiNir-misc}.

\begin{theorem}[Necessary and sufficient condition for existence and uniqueness]\label{thm:MainThm}
Let $(M,g)$ be an $n$-dimensional smooth compact Riemannian manifold with $n \geq 3$. Assume that $\Gamma$ satisfies \eqref{G1}, \eqref{G2} and that $\lambda(A_g) \in \Gamma$ in $M$. Let $S = \{p_1, \ldots, p_m\}$ be a non-empty finite subset of distinct points of $M$ and $c_1, \ldots, c_m \in (0,\infty)$.

\begin{enumerate}[(i)]
\item If $\mu_\Gamma^+ > 1$, then there exists a unique Green's function $u \in C^{0}_{loc}(M \setminus S)$ for $\Gamma$ with poles $p_1, \ldots, p_m$ and with strengths $c_1, \ldots c_m$. Furthermore, $u$ belongs to $C^{1,1}_{loc}(M \setminus S)$.
\item If $\mu_\Gamma^+ \leq 1$, Green's functions for $\Gamma$ with poles $p_i$'s and strengths $c_i$'s exist if and only if $(M,g)$ is conformal to the standard sphere and $m = 1$. 
\end{enumerate}
\end{theorem}

In Section \ref{Sec:delta}, we give a preliminary result demonstrating how Green's functions may show up in the study of blow-up sequences for nonlinear Yamabe problems.

\medskip

We list here some additional useful properties of the Green's function $u$ obtained in Theorem \ref{thm:MainThm} when $\muGp > 1$, and for given $p_i$'s and $c_i$'s.

\begin{enumerate}[(a)]
\item The Green's function $u$ is the minimum of the set of all functions $v$ in $C^{0}_{loc}(M \setminus S)$ which satisfy
\begin{align*}
&\lambda(A_{g_v}) \in \bar\Gamma \text{ and }  v > 0 \text{ in } M \setminus \{p_1, \ldots, p_m\},\\
&\lim_{x \rightarrow p_i} d_g(x,p_i)^{n-2} v(x) = c_i,\qquad i = 1, \ldots, m.
\end{align*}
See Step 2 in subsection \ref{ssec:Uniqueness}.

\item The metric $g_u$ is an asymptotically flat metric on $M \setminus \{p_1, \ldots, p_m\}$: There exists a diffeomorphism $\Phi_i$ from a punctured neighborhood of each $p_i$ into the exterior of a ball in the Euclidean space $\RR^n$ such that relative to the local coordinate functions $x^j = \Phi^j(\cdot)$ one has
\[
g_u(\partial_{x^j}, \partial_{x^\ell}) = \delta_{j\ell} + O(|x|^{-(\mu - 1)})
\]
where $\mu$ is any number in $(1,\mu_\Gamma^+] \cap (1,3)$; see Remark \ref{Rem:6.1}.

\item As a consequence of (a), Green's functions depend monotonically on $\Gamma$. More precisely, if $\Gamma \subset \Gamma'$ and $u'$ is the corresponding Green's function for $\Gamma'$ with the same poles and the same strengths, then $u \geq u'$. Similarly, the monotonicity of Green's functions with respect to the strengths $c_i$'s also holds.

\item There holds $u \geq (n-2) |\SSphere^{n-1}| \sum_{i=1}^m c_i G_{p_i}$, where $G_{p_i}$ is the Green's function for the conformal Laplacian with pole at $p_i$.

\end{enumerate}

The existence part in Theorem \ref{thm:MainThm} is proved by a suitable elliptic regularization, since equation \eqref{Eq:DegEq1} is genuinely degenerate elliptic. To solve the regularized equations as well as to show that the obtained solutions converge to a solution $u$ of \eqref{Eq:DegEq1}-\eqref{Eq:DegAs}, we construct suitable upper and lower barriers. Furthermore, our procedure allows us to construct smooth strict sub- and super-solutions of \eqref{Eq:DegEq1} which approximate the solution $u$ which we obtained. The uniqueness part then follows from a standard comparison principle argument.

\medskip
{\bf Lower Ricci bounds for continuously conformally smooth metrics.} The non-existence of smooth Green's functions when $\mu_\Gamma^+ \leq 1$ and $(M,g)$ is not conformally equivalent the standard sphere is a consequence of the rigidity of Bishop-Gromov's relative volume comparison theorem and the fact that $\lambda(A_g) \in \bar\Gamma$ with $\mu_\Gamma^+\leq 1$ implies $\Ric_g \geq 0$. In order to prove our result, we need a version of relative volume comparison theorems for continuous metrics.

When $\Gamma = \Gamma_k$ with $k > n/2$ (so that $\mu_\Gamma^+ < 1$), it was proved in the work of Gursky and Viaclovsky \cite{GV07} that Bishop-Gromov's relative volume comparison theorem (including its rigidity) holds for metrics $g_{u} \in C^{1,1}_{\rm loc}$ where $u$ is the $C^{1,\alpha}_{loc}$ limit of a sequence of smooth functions $u_j$ which are bounded in $C^2_{loc}$ and satisfy $\lambda(A_{g_{u_j}}) \in \bar\Gamma$. 

Our treatment for Bishop-Gromov's relative volume comparison theorem is different from \cite{GV07}. Note that our definition of Green's functions $u$ only gives the continuity of the metric $g_u$. We exploit the fact that metrics of interest to us are conformal to smooth metrics, which we will refer to as continuously conformally smooth metrics. For this class of metrics, we can define a notion of (purely local) lower Ricci curvature bounds in the sense of viscosity; see Definition \ref{Def:LRBVis}. This is naturally coherent with the notion of viscosity (super-)solutions for \eqref{Eq:DegEq1}. We establish the following purely local relative volume comparison theorem (see Section \ref{Sec:VolComp} for terminologies):

\begin{theorem}[Relative volume comparison]\label{Thm:RelVolComp}
Let $(M^n,g)$ be a smooth complete Riemannian manifold of dimension $n \geq 2$, $f \in C^0_{loc}(M)$, and $k$ be a constant. Suppose $\Ric(e^{2f}g) \geq (n-1)k$ in some ball $B_{e^{2f}g}(p,R)$ centered at $p$ and of radius $R > 0$ with respect to the metric $e^{2f}g$ in the viscosity sense. If $k > 0$, suppose further that $R \leq \frac{\pi}{2\sqrt{k}}$. Then, for $r \in (0,R)$,  the function
\[
r \mapsto \frac{Vol_{e^{2f}g}(B_{e^{2f}g}(p,r))}{v(n,k,r)}
\]
is a non-increasing function, where $v(n,k,r)$ is the volume of a ball of radius $r$ in the simply connected constant curvature space form $\SSphere^n_k$. 

In addition, if it holds for some $p \in M$ and $r > 0$ (and $8r \leq \frac{\pi}{2\sqrt{k}}$ if $k > 0$) that $Vol_{e^{2f}g}(B_{e^{2f}g}(p,8r)) = v(n,k,8r)$, then $f$ is smooth in $B_{e^{2f}g}(p,r)$ and $B_{e^{2f}g}(p,r)$ is isometric to a ball of radius $r$ in the simply connected constant curvature space form $\SSphere^n_k$.
\end{theorem}

It would be interesting to relate our notion of lower Ricci bounds in the viscosity sense to notions of lower Ricci bounds related to Bakry-\'Emery inequalities or convexity of entropies. For the latter, see Ambrosio \cite{AmbrosioICM18} and the references therein.

\medskip
{\bf Asymptotics of Green's functions.} It is well known, in the case of the scalar curvature, that the Green's function $G_p$ can arise as the limit of a suitably rescaled blow-up sequence of solutions to the Yamabe problem. This limit object $G_p$ has an asymptotic expansion near $p$ (cf. Lee and Parker \cite{Lee-Parker}) which contains local as well as global geometric information about $(M,g)$. In particular, in a conformal normal coordinate system at $p$, when $3 \leq n \leq 5$ or when the Weyl tensor of $g$ vanishes suitably fast near $p$, we have
\[
G_p(x) 
	= \frac{1}{(n-2) |\SSphere^{n-1}|} \big(d_g(x,p)^{-(n-2)}  + A_p + O(d_g(x,p))\big).
\]
In such case, the metric $G_p^{\frac{4}{n-2}}g$ is asymptotically flat and scalar flat on $M \setminus\{p\}$ and its ADM mass is, up to a dimensional constant, the constant $A_p$ in the above expansion. The positivity of the ADM mass plays an important role in the resolution of the compactness problem for the Yamabe problem (see Brendle and Marques \cite{BrendleM09}, Khuri, Marques and Schoen \cite{KMS09} and the references therein) and more generally in the study of scalar curvature.

It is therefore of interest to study Green's functions and their asymptotic behaviors in the current fully nonlinear setting, and, in particular, to understand what geometric information they encode. The following result gives a first step in this direction. Since its proof is of different nature than what is being discussed in this paper, it will appear elsewhere.

\begin{theorem}[Estimates for Green's functions]\label{Thm:21II19-EstGF}
Let $(M,g)$ be an $n$-dimensional smooth compact Riemannian manifold with $n \geq 3$. Assume that $\Gamma$ satisfies \eqref{G1}, \eqref{G2}, $\mu_\Gamma^+ > 1$ and that $\lambda(A_g) \in \Gamma$ in $M$. Let $S = \{p_1, \ldots, p_m\}$ be a non-empty finite subset of distinct points of $M$ and $c_1, \ldots, c_m \in (0,\infty)$, and $u \in C^{0}_{\rm loc}(M \setminus S)$ be the Green's function for $\Gamma$ with poles $p_1, \ldots, p_m$ and with strengths $c_1, \ldots, c_m$. Then $u \in C^{1,1}_{\rm loc}(M \setminus S)$ and there exist constants $\kappa > 0, r_0 > 0$ and $C > 0$ such that, for $i = 1, \ldots, m$ and $x \in B(p_i, r_0)$, there hold
\begin{align}
|u(x) - c_i d_g(x,p_i)^{2-n}| &\leq Cd_g(x,p_i)^{2-n + \kappa},\label{Eq:28I19-A1}\\
|\nabla(u(x) - c_i d_g(x,p_i)^{2-n})| &\leq Cd_g(x,p_i)^{1-n + \kappa},\label{Eq:28I19-A2}
\end{align}
and
\begin{align}
|\nabla^2 u(x)| \leq Cd_g(x,p_i)^{-n} .\label{Eq:28I19-A3}
\end{align}
If it holds in addition that $(1, 0, \ldots, 0) \in \Gamma$, then
\begin{align}
|\nabla^2(u(x) - c_i d_g(x,p_i)^{2-n})| &\leq Cd_g(x,p_i)^{-n + \kappa}.\label{Eq:28I19-A4}
\end{align}
\end{theorem}

It would be interesting to see if estimate \eqref{Eq:28I19-A4} holds for all cones $\Gamma$ (with $1 < \mu_\Gamma^+ <  n-1$), or at least for $\Gamma_k$ with $2 \leq k < \frac{n}{2}$. It is readily seen that the metric $g_u$ is asymptotically flat. If estimates \eqref{Eq:28I19-A1}-\eqref{Eq:28I19-A2} and \eqref{Eq:28I19-A4} hold for $\Gamma = \Gamma_k$ for some $\kappa > \frac{n-2k}{k+1}$, then it can be shown that $g_u$ has a well-defined $k$-mass (see Li and Nguyen \cite{LiNgkMass} and Ge, Wang and Wu \cite{GeWangWu-AdvM14}). It is of much interest to study whether a generalized mass (as in \cite{GeWangWu-AdvM14, LiNgkMass}, or a variant of such) can be defined for $g_u$ (including the case $\Gamma = \Gamma_k$), what role it plays, or whether it enjoys a similar positive mass result, etc. (Note that, when $(M,g)$ is locally conformally flat and not conformally equivalent to the standard sphere, the positivity of mass is a consequence of \cite[Theorem 1.2]{LiNgBocher}. See also \cite{GeWangWu14-IMRN} when $(M,g)$ is conformally flat.)

In a sense, the gradient and Hessian estimates in Theorem \ref{Thm:21II19-EstGF} can be viewed as ones for `the linearized equation' of \eqref{Eq:DegEq1} near the fundamental solution. We believe that these estimates and their analogues for solutions to \eqref{Eq:12VIII20-E1} which are close to the fundamental solution, if hold, should be of importance in understanding compactness issues for \eqref{Eq:12VIII20-E1}.

Similar notions of Green's functions for fully nonlinear elliptic Hessian-type equations have been studied in the literature; see e.g. Armstrong, Sirakov and Smart \cite{ArmsSirSmart11}, Harvey and Lawson  \cite{HarveyLawson16-AdvM}, Jin and Xiong \cite{JinXiong16-AdvM}, J\"orgens \cite{Jorgens55-MAnn}, Labutin \cite{Labutin01}, Trudinger and Wang \cite{TWHessianII}. We mention here a recent paper by Esposito and Malchiodi \cite{EspositoMalchiodi19} where a related result was established in a context involving log-determinant functionals.

We conclude the introduction with the following question:
\begin{question}
Is the $C^{1,1}_{\rm loc}$ viscosity solution to \eqref{Eq:DegEq1}-\eqref{Eq:DegAs} constructed in Theorem \ref{thm:MainThm}  smooth in a punctured neighborhood of the $p_i$'s, at least for $\Gamma = \Gamma_2$? 
\end{question}
This question is motivated by a result of Lempert \cite{Lempert83-MathAnn}, which asserts that for any strictly convex and analytically bounded $\Omega \subset \mathbb{C}^n$, any real analytic $\varphi: \partial \Omega \rightarrow \mathbb{R}$, and any $p_0 \in \Omega$, there exists $C_0 > 0$ such that for all $C > C_0$ there exists a unique solution, real analytic  in $\Omega \setminus \{p_0\}$ and pluri-subharmonic in $\Omega$, to the degenerate complex Monge-Amp\`ere problem $(\partial \bar\partial)^n u = 0$ in $\Omega \setminus \{p_0\}$, $u(z) = C\ln |z - p_0| + O(1)$ and $u = \varphi$ on $\partial \Omega$.

The rest of the paper is structured as follows. In Section \ref{Sec:VolComp}, we define a suitable notion of lower Ricci curvature bounds for continuous metrics which are conformal to smooth metrics and prove a version of Bishop-Gromov's relative volume comparison theorem for these metrics. In Section \ref{Sec:NoExistence}, we use the relative volume comparison theorem to prove part (ii) of Theorem \ref{thm:MainThm}. The proof of part (i) of Theorem \ref{thm:MainThm} is then carried out in Section \ref{Sec:ExistenceUniq}. Section \ref{Sec:delta} is devoted to a result illustrating the relevance of Green's functions in the study of nonlinear Yamabe problems. The paper includes also two appendices, one on the construction of a concave function whose zeroth level set is $\partial\Gamma$ and another one on the convexity of the set of eigenvalues of matrices belonging to a convex set.

\section{Lower Ricci bounds for continuously conformally smooth metrics and volume comparison}\label{Sec:VolComp}

In this section, we introduce a notion of lower Ricci bounds in the viscosity sense for continuous metrics which are conformal to smooth metrics. As mentioned in the introduction, we will conveniently refer to these metrics as continuously conformally smooth metrics. We establish Theorem \ref{Thm:RelVolComp}, a version of Bishop-Gromov's relative comparison theorem. This will be used to prove statement (ii) in Theorem \ref{thm:MainThm}, i.e. the non-existence of solutions to \eqref{Eq:DegEq1}-\eqref{Eq:DegAs} when $\muGp \leq 1$.

It is instructive to note the fact that if $g$ is a smooth metric and $f$ is a smooth function, then a bound for the Ricci tensor of the conformal metric $e^{2f}g$ translates to a second order, though nonlinear, partial differential inequality for the function $f$. One can thus define the notion of a lower Ricci curvature bound for $e^{2f}g$ when $f$ is merely continuous in the viscosity sense, as one does for nonlinear second order elliptic equations. See Definition \ref{Def:LRBVis}.

A nice feature of this way of defining lower Ricci curvature bounds is that if a metric has a lower Ricci curvature bound, then it can be approximated by locally Lipchitz metrics which also satisfy related lower Ricci curvature bounds. See Proposition \ref{Prop:InfConvStab}.

We then proceed to approximate locally Lipschitz metrics with lower Ricci curvature bounds by smooth metrics. While it is desirable to keep a pointwise lower Ricci curvature bound for the approximants, we are content with keeping a suitable integral lower Ricci curvature bound. See Proposition \ref{Prop:LipWeakStab}. The relative volume comparison is then drawn from results of Peterson and Wei \cite{PetersenWei-GAFA97, PetersenWei-TrAMS01}, Wei \cite{Wei-SDG07} on smooth metrics of integral lower Ricci curvature bounds.

Last but not least, a subtle point in the proof of the rigidity of relative volume comparison is to prove that the metric-space isometry between the given continuous metric and the corresponding constant curvature metric is a smooth Riemannian isometry with respect to the given smooth structure. We again exploit the fact that the given continuous metric is conformal to a smooth metric and show that the isometry in fact satisfies the $n$-Laplacian equation, which is the Euler-Lagrange equation of a conformally invariant functional. We then appeal to the regularity theory for the $n$-Laplacian to reach the conclusion.

\subsection{Two notions of Ricci lower bounds}

Assume that $g$ is a smooth metric on a smooth (compact or non-compact) manifold $M^n$ of dimension $n \geq 2$ and $f$ is a continuous function defined on an open subset $\Omega \subseteq  M$. Let us first start by defining what we mean by a lower Ricci bound for $e^{2f}g$. 

\begin{definition}\label{Def:LRBVis}
Let $k$ and $f$ be continuous functions defined on an open subset $\Omega$ of a smooth Riemannian manifold $(M,g)$. We say that $\Ric(e^{2f}g) \geq (n-1)k$ in the viscosity sense in $\Omega$ if, for every $x_0 \in \Omega$ and for every $\varphi \in C^2(\Omega)$ such that $\varphi(x_0) = f(x_0)$ and $\varphi \leq f$ in a neighborhood of $x_0$, one has 
$$\Ric(e^{2\varphi}g)(x_0) - (n-1)k(x_0)\,e^{\varphi(x_0)}g(x_0) \text{ is non-negative definite}.$$
\end{definition}

It is clear that $\Ric(e^{2f}g) \geq (n-1)k$ in the viscosity sense if and only if it holds for any continuous non-negative definite $(2,0)$-tensor $a$ defined on $\Omega$ that
\begin{align*}
a^{ij}\Ric_{ij}(e^{2f}g) 
	&= -(n-2) a^{ij}\nabla_{ij} f - \tr_g(a)\Delta_g f\\
		&\qquad\qquad  + (n-2) a(df,df) - (n-2)|df|_g^2\,\tr_g(a) + a^{ij}\Ric_{ij}(g) \\
	&\geq (n-1)k \,\tr_g(a)
\end{align*}
in the usual viscosity sense. In addition, if $f$ is $C^{2}$ and satisfies $\Ric(e^{2f}g) \geq (n-1)k$ in the viscosity sense, then $\Ric(e^{2f}g) \geq (n-1)k$ in the classical sense.

If $f$ is Lipschitz continuous, the quadratic term in the expression for $\Ric(e^{2f}g)$ is integrable. This motivates the following definition.

\begin{definition}
Let $k$ be a continuous function and $f$ be a locally Lipschitz function defined on an open subset $\Omega$ of a smooth Riemannian manifold $(M,g)$. We say that $\Ric(e^{2f}g) \geq (n-1)k$ in the weak sense in $\Omega$ if, for every smooth compactly supported non-negative definite $(2,0)$-tensor $a$ defined on $\Omega$, there holds
\begin{eqnarray}
&&\int_\Omega \Big[(n-2) \nabla_{i} f \,\nabla_j a^{ij} + \nabla_i f\,\nabla^i \tr_g(a)
	+ (n-2) a(df,df) - (n-2)|df|_g^2\,\tr_g(a) \nonumber\\
	&&\qquad\qquad  + a^{ij}\Ric_{ij}(g) \Big]\,dv_g
		\geq \int_\Omega (n-1)k \,\tr_g(a)\,dv_g .\label{Eq:wRicLB}
\end{eqnarray}
\end{definition}

We will prove later that, if $f$ is Lipschitz and if $\Ric(e^{2f}g) \geq (n-1)k$ in the viscosity sense, then $\Ric(e^{2f}g) \geq (n-1)k$ in the weak sense; see Proposition \ref{Prop:Vis=>Weak}.

One key property concerning metrics with lower Ricci bounds in the viscosity sense which we will establish is the following result. Roughly speaking, every continuous metric $e^{2f}g$ whose Ricci curvature is bounded from below in the viscosity sense can be approximated by smooth conformal metrics $e^{2\bar f_\ell}g$ whose Ricci curvatures are bounded from below in $L^p$-sense for all $1 \leq p < \infty$. More precisely, we prove:

\begin{proposition}\label{Prop:Approx}
Let $\Omega$ be a bounded open subset of $M$ and $f, k \in C(\bar\Omega)$ such that $\Ric(e^{2f}g) \geq (n-1)k$ in the viscosity sense in $\Omega$. Then there exists a sequence of functions $\{\bar f_\ell\} \subset C^\infty(\Omega)$ which converges locally uniformly to $f$ such that, for any $1 \leq p < \infty$ and open $\omega \Subset \Omega$, 
\[
\lim_{\ell \rightarrow \infty} \int_{\omega} \Big\{ \max\big(-\lambda_1(\Ric(e^{2\bar f_\ell}g)) + (n-1)k, 0\big)\Big\}^p\,dv_g = 0,
\]
where $\lambda_1(\Ric(e^{2\bar f_\ell}g))$ is the smallest eigenvalue of $\Ric(e^{2\bar f_\ell}g)$ with respect to $e^{2\bar f_\ell}g$.
\end{proposition}

\begin{proof} This is an immediate consequence of Propositions \ref{Prop:InfConvStab} and \ref{Prop:LipWeakStab} below concerning the stability of our notion of Ricci lower bounds under two different regularization processes: the inf-convolution and the convolution against a kernel.
\end{proof}

\subsection{Stability of Ricci lower bounds under inf-convolutions}

In this section, we prove that every continuous metric $e^{2f}g$ whose Ricci curvature is bounded from below in the viscosity sense can be approximated by Lipschitz conformal metrics $e^{2\bar f_\ell}g$ whose Ricci curvatures are also bounded from below in the viscosity sense. We prove:

\begin{proposition}\label{Prop:InfConvStab}
Let $(M,g)$ be a smooth complete Riemannian manifold. Let $\Omega$ be a bounded open subset of $M$ and $f,k \in C(\bar\Omega)$ such that $\Ric(e^{2f}g) \geq (n-1)k$ in the viscosity sense in $\Omega$. Then, for all sufficiently small $\vareps > 0$, there exist functions $f_\vareps \in C^{0,1}_{\rm loc}(\Omega) \cap C(\bar\Omega)$ which are locally semi-concave and $\theta_\vareps \in C(\Omega)$ such that $\Ric(e^{2f_\vareps}g) \geq (n-1)k - \theta_\vareps$ in $\Omega$ in the viscosity sense, $f_\vareps \rightarrow f$ in $C(\bar\Omega)$ and $\theta_\vareps \rightarrow 0$ in $C^0_{loc}(\Omega)$ as $\vareps \rightarrow 0$.
\end{proposition}

We will use inf-convolutions to regularize. Let $\Omega$ be a bounded open subset of $M$. For $f \in C(\bar\Omega)$ and small $\vareps > 0$, we define 
\[
f_\vareps(x) = \inf_{y \in \Omega} \Big[f(y) + \frac{1}{\vareps}d_g(x,y)^2\Big], \qquad x \in \Omega,
\]
where $d_g$ denotes the distance function of $g$. We note that $f_\vareps$ satisfies the following properties; see e.g. \cite[Chapter 5]{CabreCaffBook} or \cite[Section 2]{LiNgWang}.

\begin{enumerate}[(i)]
\item\label{UpLowPropi} $f_\vareps \in C(\bar\Omega)$ is monotonic in $\vareps$ and $f_\vareps \rightarrow f$ uniformly as $\vareps \rightarrow 0$.

\item\label{UpLowPropii} $f_{\vareps}$ is punctually second order differentiable (see e.g. \cite{CabreCaffBook} for a definition) almost everywhere in $\Omega$ and $\nabla^{2}f_{\vareps}\leq C(\Omega,g) \vareps^{-1} g$ a.e. in $\Omega$.

\item\label{UpLowPropiii} For any $x \in \Omega$, there exists $x_* = x_*(x) \in \bar \Omega$ such that 
\begin{equation}
f_\vareps(x) = f(x_*)  + \frac{1}{\vareps}d_g(x,x_*)^2 .\label{Eq:UpLowx*Def}
 \end{equation}

\item \label{UpLowPropiv} For any non-empty open subset $\omega$ of $\Omega$, there holds
\[
|\nabla f_\vareps| \leq \frac{C(\Omega,g)}{\vareps^{\frac{1}{2}}} \big[\sup_{\omega} f - \min_{\bar\Omega} f\big]^{\frac{1}{2}} \text{ a.e. in } \omega.
\]

\item \label{UpLowPropv} If $|f(x) - f(y)| \leq m(d_g(x,y))$ for all $x, y \in \bar\Omega$ and for some non-negative continuous non-decreasing function $m: [0,\infty) \rightarrow [0,\infty)$ satisfying $m(0) = 0$, then
\begin{equation}
d_g(x,x_*) \leq \big[\vareps\,m((C(\Omega,g)\vareps\,\sup_{\bar\Omega} |f|)^{1/2})\big]^{1/2}.
\label{Eq:UpLowdxx*}
\end{equation}
\end{enumerate}

\begin{proof}[Proof of Proposition \ref{Prop:InfConvStab}] Since $\bar\Omega$ is compact, it is enough to consider the case that $\Omega$ is contained in a single chart of $M$. Fix a compact subset $\omega$ of $\Omega$ and a point $x^0 \in \omega$. We will prove that for every $\varphi \in C^2(\Omega)$ such that $\varphi\leq f_\vareps$ in a neighborhood of $x^0$ and $\varphi(x^0) = f_\vareps(x^0)$ it holds 
\begin{equation}
\Ric(e^{2\varphi}g)(x^0) \geq [(n-1)k(x^0) - o_\vareps(1)]\,e^{2\varphi}(x^0)\,g(x^0),
	\label{Eq:12VII18-E1}
\end{equation}
where here and below $o_\vareps(1)$ denotes some constant which depends only on $\vareps$, $\|f\|_{C(\bar\Omega)}$, $dist(\omega,\partial\Omega)$ and the moduli of continuity of $f$ and $k$ on $\bar\Omega$ such that $o_\vareps(1) \rightarrow 0$ as $\vareps \rightarrow 0$.

By the definition of $f_\vareps$, $f(x) \geq f_\vareps(y) - \frac{1}{\vareps} d_g(x,y)^2$ for all $x, y \in \Omega$. Thus, for $x, y$ close to $x^0$,
\[
f(x) \geq \varphi(y) - \frac{1}{\vareps} d_g(x,y)^2.
\]
Now if $x^0_* = x_*(x^0)$ is defined as in  \eqref{Eq:UpLowx*Def} and if $y$ is a $C^2$ map defined from on a neighborhood of $x^0_*$ into $\Omega$ such that $y(x^0_*) = x^0$, then 
\[
f(x) \geq \varphi(y(x)) - \frac{1}{\vareps} d_g(x,y(x))^2 =: \psi(x)\text{ near } x^0 \text{ and } f(x^0) = \psi(x^0).
\]
Hence, as $\Ric(e^{2f}g) \geq (n-1)k$ in the viscosity sense, we have that
\begin{equation}
\Ric(e^{2\psi}g)(x^{0}_*)  \geq (n-1)k(x^{0}_*) e^{2\psi(x^{0}_*)} g.
	\label{Eq:13VII18-A1}
\end{equation}
We will deduce \eqref{Eq:12VII18-E1} from \eqref{Eq:13VII18-A1} by a judicious choice of $y$.

For expository purpose and to motivate our later argument, let us first present the case where $\Omega$ is a Euclidean domain and $g$ is the Euclidean metric. The general case will be treated subsequently. 

When $g$ is the Euclidean metric, 
\[
\Ric(e^{2\varphi}g) = -(n-2)\nabla^2 \varphi - \Delta \varphi\,Id + (n-2)d\varphi \otimes d \varphi - (n-2)|d \varphi|^2\,Id.
\]
Now let
\[
y(x) = x + x^0 - x^0_*
\]
so that $\nabla \psi(x^{0}_*) = \nabla \varphi(x^{0})$ and $\nabla^2 \psi(x^{0}_*) = \nabla^2 \varphi(x^{0})$. Estimate \eqref{Eq:12VII18-E1} is therefore readily seen from  \eqref{Eq:13VII18-A1} and \eqref{Eq:UpLowdxx*}.

Let us now turn to the case when $g$ is a general Riemannian metric. The proof above uses strongly the fact that, when $(\Omega,g)$ is Euclidean, the tangent and cotangent spaces of $M$ at $x^{0}$ and $x^{0}_*$ can be naturally identified and this identification does not interfere with the equation. This has the advantage that in our choice of the function $y$, the $\vareps$-dependent contribution in the test function $\psi$ is a constant. In the general setting, special care must be given.

An inspection leads to the following choice of $y$:
\[
y(\exp_{x^{0}_*}(z_*)) = \exp_{x^{0}}(z)
\]
where $z = Pz_* \in T_{x^{0}} M$ and $P: T_{x^{0}_*}M \rightarrow T_{x^{0}}M$ is the parallel transport map along the (unique) minimizing geodesic $\gamma_{x^0_*, x^0}$ connecting $x^{0}_*$ to $x^{0}$. The map $y$ translates a neighborhood of $x^0_*$ to that of $x^0$ along the geodesic $\gamma_{x^0_*, x^0}$.

By the first and second variation formulae for length (see e.g. \cite[Theorems 3.31 and 3.34]{GalHulLaf}), we have that
\begin{align*}
\frac{d}{dt}\Big|_{t = 0} d_g(\exp_{x^{0}}(t z), \exp_{x^{0}_*}(t z_*)) 
	& = 0,\\
\frac{d^2}{dt^2} d_g(\exp_{x^{0}}(t z), \exp_{x^{0}_*}(t z_*)) 
	& 
	   = 	O((d_g(x^{0}, x^{0}_*) + |t||z_*|_g) |z_*|_g^2) \text{ for small $|t|$}.
\end{align*}
Hence
\[
d_g(\exp_{x^{0}}(z), \exp_{x^{0}_*}( z_*))
	= d_g(x^{0}, x^{0}_*) + O(d_g(x^{0}, x^{0}_*))|z_*|_g^2,
\]
and so
\[
\psi(\exp_{x^{0}_*}(z_*)) = \varphi(\exp_{x^{0}}(z)) -  \frac{1}{\vareps} d_g(x^{0} ,x^{0}_*)^2  + o(|z_*|_g^2).
\]
Loosely speaking, this means that the $\vareps$-dependent contribution in the test function $\psi$ is constant up to a super-quadratic error. (In fact, the choice of $y$ which ensures this property is unique up to quadratic terms in the Taylor expansion of $y$ around $x^0_*$.) We hence obtain
\begin{align}
d\psi(x^{0}_*)(z_*) 
	&= d \varphi(x^{0})(z),\label{Eq:14VII18-B1} \\
\nabla_g^2 \psi(x^{0}_*)(z_*,z_*) 
	&= \nabla_g^2 \varphi(x^{0})(z,z). \label{Eq:14VII18-B2}
\end{align}

Now, recall from \eqref{Eq:13VII18-A1} that
$$\Ric(e^{2\psi}g)\big|_{x^{0}_*}(z_*,z_*) \geq (n-1)k(x^{0}_*) e^{2\psi(x^{0}_*)}|z_*|_g^2.$$
Using \eqref{Eq:14VII18-B1}, \eqref{Eq:14VII18-B2} and the fact that the transformation $z_* \mapsto Pz_* = z$ (from $T_{x^{0}_*}M$ to $T_{x^{0}}M$) is length preserving, we obtain
\begin{align*}
 &-(n-2)\nabla_g^2 \varphi(x^{0})(z,z)  - \Delta_g \varphi(x^{0})\,|z|_g^2\\
 	&\qquad + (n-2)[d\varphi(x^0) (z)]^2 - (n-2)|d\varphi|_g^2(x^{0})\,|z|_g^2 + \Ric(g)\big|_{x^{0}_*}(z_*,z_*) \\
	&\qquad\qquad \geq (n-1)k(x^{0}_*)\,e^{2\varphi(x^{0}) - \frac{2}{\vareps} d_g(x^{0} ,x^{0}_*)^2} |z|_g^2.
\end{align*}
Recalling \eqref{Eq:UpLowdxx*}, we obtain \eqref{Eq:12VII18-E1}, which concludes the proof.
\end{proof}

\subsection{Viscosity Ricci lower bounds imply weak Ricci lower bounds for Lipschitz conformal factors}

In this subsection, we prove:
\begin{proposition}\label{Prop:Vis=>Weak}
Let $(M,g)$ and $\Omega$ be as in Proposition \ref{Prop:InfConvStab}. Assume that $f \in C^{0,1}_{loc}(\Omega)$ and $k \in C(\bar\Omega)$. If $\Ric(e^{2f}g) \geq (n-1)k$ holds in the viscosity sense in $\Omega$, then it holds in the weak sense.
\end{proposition}

\begin{proof} Without loss of generality, we can assume that $\Omega$ is bounded, $\partial\Omega$ is smooth, $k \in C^0(\bar\Omega)$, and $f \in C^{0,1}(\bar\Omega)$. Furthermore, by using Proposition \ref{Prop:InfConvStab}, we may further assume that $f$ is almost everywhere punctually second order differentiable and that $\nabla^2 f \leq C$ a.e. in $\Omega$.

We will establish \eqref{Eq:wRicLB} for an arbitrary smooth $(2,0)$-tensor $a$ defined on $\bar\Omega$ such that $a \equiv 0$ on $\partial\Omega$. Writing $a = \sum_k \psi_k a$ for a suitable partition of unity $\{\psi_k\}$ if necessary, it suffices to consider the case that $\Omega$ is contained in a single chart. Furthermore, by considering $a + \delta \varphi g^{-1}$ (instead of $a$) for all sufficiently small $\delta > 0$ and some $\varphi\in C^\infty(\bar\Omega)$ satisfying $\varphi > 0$ in $\Omega$ and $\varphi = 0$ on $\partial\Omega$, we may assume that $a$ is positive definite in $\Omega$.

Set $b^{ij} = (n-2)a^{ij} + \tr_g(a)\,g^{ij}$ and
\begin{multline*}
h = \nabla_j b^{ij}\,\nabla_i f - (n-2) a(df,df) + (n-2)|df|_g^2\,\tr_g(a)\\
		  - a^{ij}\Ric_{ij}(g) 
		+  (n-1)k \,e^{2f} \,\tr_g(a) \in L^\infty(\Omega).
\end{multline*}
We note that $(b^{ij})$ is positive definite in $\Omega$. Since the subdifferential map of a convex function has a closed graph (see e.g. \cite[Theorem 24.4]{Rockafellar}) and is single-valued almost everywhere in its domain, we can, without loss of generality, identify $h$ with its lower semi-continuous representative. 

To prove \eqref{Eq:wRicLB}, we show that 
\begin{equation}
-\nabla_i (b^{ij} \nabla_j f) \geq h \text{ in } \Omega \text{ in the weak sense.}
	\label{Eq:05I18-X1}
\end{equation}

\medskip
\noindent{\it Step 1:} We start with showing a comparison principle for $f$. For a subdomain $\omega \subset \Omega$ with smooth boundary $\partial\omega$, let $v_\omega$ be the solution to
\begin{align*}
L(v_\omega):= - \nabla_i (b^{ij} \nabla_j  v_\omega)
	&=  h \text{ in } \omega
\end{align*}
subjected to the Dirichlet boundary condition $v_\omega = f$ on $\partial\omega$. We claim that $v_\omega \leq f$ in $\omega$. 

Indeed, since $h$ is lower semi-continuous, there exists a sequence of smooth functions $h_l \leq h - \frac{1}{l}$ which converges pointwise to $h$ as $l \rightarrow \infty$. Let $v_l$ solves
\begin{align*}
L(v_l) 
	&=  h_l \text{ in } \omega,\\
v_l
	&= f - \frac{1}{l} \text{ on } \partial\omega.
\end{align*}
To prove the claim it suffices to show that $m_l := \inf_\omega (f - v_l) \geq 0$. Assume by contradiction that $m_l < 0$. Pick some small $\eta \in (0,|m_l|)$ for the moment and let $\xi = \xi_{l,\eta}= f - v_l - m_l - \eta$ and $\Gamma_\xi$ be the convex envelope of $-\xi^- = -\max(-\xi,0)$. By the Alexandrov-Bakelman-Pucci estimate \cite[Lemma 3.5]{CabreCaffBook} (which applies since $f$ is semi-concave and $\Omega$ is contained in a single chart), the set $\{\xi = \Gamma_\xi\}$ has non-empty measure. Thus there is a point $x_{l,\eta}$ in this set where $\xi$ is punctually second order differentiable and
\begin{equation}
-\eta \leq \xi(x_{l,\eta}) \leq 0, |\partial \xi(x_{l,\eta})| \leq C \eta, \text{ and } \partial^2 \xi(x_{l,\eta}) \geq 0,
	\label{Eq:03I18-X1}
\end{equation}
where $\partial$ denotes the partial derivatives and $C$ is independent of $\eta$. At this point, $f$ is punctually second order differentiable and so
\[
\Ric(e^{2f}g)(x_{l,\eta}) \geq (n-1)k(x_{l,\eta})e^{2(f(x_{l,\eta}) + m_l)}g(x_{l,\eta}),
\]
which implies
\[
L f(x_{l,\eta}) \geq h(x_{l,\eta}) \geq h_l(x_{l,\eta}) + \frac{1}{l}.
\]
In view of \eqref{Eq:03I18-X1}, this implies that
\[
Lv_l(x_{l,\eta}) \geq h_l(x_{l,\eta}) + \frac{1}{l} - C\eta > h_l(x_{l,\eta}),
\]
provided $\eta$ is chosen sufficiently small. This contradicts the definition of $v_l$. The claim is proved.

\medskip
\noindent{\it Step 2:} We now proceed to prove \eqref{Eq:05I18-X1}. Fix a sequence of smooth functions $\{f_l\} \subset C^\infty(\bar\Omega)$ which converges uniformly to $f$ in $\bar\Omega$ and satisfies $f_l < f$ in $\bar\Omega$. Fix some subdomain $\omega \Subset \Omega$ with smooth boundary $\partial\omega$. Let $\xi_l$ be the solution to the (obstacle) variational problem
\[
\min \Big\{\int_\omega [b^{ij} \nabla_i \xi \nabla_j \xi  - h\xi]\,dv_g: \xi \in H^1(\omega), \xi\big|_{\partial\omega} = f_l\big|_{\partial\omega}, \xi \geq f_l \text{ in } \omega\Big\}.
\]
It is well known that the minimizer $\xi_l$ to the above problem exists uniquely and $\xi_l$ satisfies
\[
L(\xi_l) \geq h \text{ in the weak sense in $\omega$},
\]
and
\[
L(\xi_l) = h \text{ in the weak sense in $\{\xi_l > f_l\}$}.
\]
Hence, by Step 1, we have 
\[
\xi_l \leq f \text{ in } \omega.
\]
Consequently, by the uniform convergence of $f_l$ to $f$, we have that $\{\xi_l\}$ converges uniformly to $\xi$ on $\omega$ and so $Lf \geq h$ in the weak sense in $\omega$. Since $\omega$ is arbitrary, we have thus proved \eqref{Eq:05I18-X1}.
\end{proof}

\subsection{Stability of Ricci lower bounds under convolutions against a smooth kernel}

We have seen above that the inf-convolution `preserves' Ricci lower bounds and improves the regularity of conformal factors from continuity to Lipschitz continuity. In this subsection, we are concerned with approximations with better regularity.

Throughout this subsection, we assume that $f \in C^{0,1}_{loc}(\Omega)$ unless otherwise stated.

Let $\varrho: \RR \rightarrow [0,\infty)$ be an even smooth function of compact support such that 
$$|\SSphere^{n-1}|\int_0^\infty t^{n-1}\,\varrho(t)\,dt = 1,$$
and define $\varrho_\vareps(t) = \vareps^{-n} \varrho(\vareps^{-1}t)$. A smoothing $\{\bar f_\vareps\}$ of $f$ is then obtained by convolution against $\varrho_\vareps$:
\[
\bar f_\vareps(x) = \int_{M} \varrho_\vareps(d(x,y))\,f(y)\,dv_g(y) \text{ for } x \in \Omega_\vareps := \{x \in \Omega: d(x,\partial\Omega) > \vareps\}.
\]
Noting that
\begin{equation}
Z_\vareps(x) := \int_{M} \varrho_\vareps(d(x,y))\,dv_g(y) \rightarrow 1 \text{ in } C^2_{loc}(M),
	\label{Eq:05I18-ZepsDef}
\end{equation}
we see that $\bar f_\vareps \rightarrow f$ in $C^{0,\alpha}_{loc}(\Omega)$ for any $\alpha \in (0,1)$ and $\nabla \bar f_\vareps \rightarrow \nabla f$ a.e. in $\Omega$.

The following result establishes the stability of pointwise Ricci lower bounds for $C^{0,1}$ conformal metrics. For Lipschitz conformal metrics, we prove an integral stability statement, which suffices for our purpose.
 
\begin{proposition}\label{Prop:LipWeakStab}
Let $(M,g)$ and $\Omega$ be as in Proposition \ref{Prop:InfConvStab}. Assume that $f \in C^{0,1}_{loc}(\Omega)$, $k \in C(\bar\Omega)$ and $\Ric(e^{2f}g) \geq (n-1)k$ in the viscosity sense in $\Omega$. Then, for any $1 \leq p < \infty$ and open $\omega \Subset \Omega$, the smallest eigenvalue $\lambda_1(\Ric(e^{2\bar f_\vareps}g))$ of $\Ric(e^{2\bar f_\vareps}g)$ with respect to $e^{2\bar f_\vareps}g$ satisfies
\[
\lim_{\vareps \rightarrow 0} \int_{\omega} \Big\{ \max\big(-\lambda_1(\Ric(e^{2\bar f_\vareps}g)) + (n-1)k, 0\big)\Big\}^p\,dv_g = 0.
\]
In addition, if $f \in C^1(\bar\Omega)$, then, for all sufficiently small $\vareps > 0$, there exists $\theta(\vareps) > 0$ (which possibly depends on $\omega$, $f$ and $k$) such that $\Ric(e^{2\bar f_\vareps}g) \geq (n-1)k - \theta(\vareps)$ in $\omega$ and $\theta(\vareps) \rightarrow 0$ as $\vareps \rightarrow 0$. 
\end{proposition}

The very rough idea of the proof is as follows. Ignoring lower derivatives, one can roughly think of a lower bound for $\Ric(e^{2f}g)$ as a requirement that the Hessian of $f$ belongs to certain convex subset in the bundle of symmetric $(0,2)$-tensors. The convolution is in fact an averaging process and thus, in principle, preserves such convexity. For example, Greene and Wu showed in \cite[Proposition 2.2]{GreeneWu79-AnnENS} that continuous geodesically strictly convex functions can be approximated by smooth geodesically strictly convex functions. As we are dealing with `convexity constraint' in the viscosity sense, the argument in \cite{GreeneWu79-AnnENS} does not apply directly. In fact our proof below does not work if we relax $f \in C^{0,1}$ to $f \in C^0$.

Before establishing a lower Ricci bound for the metric $e^{2\bar f_\vareps}g$, we briefly discuss some facts about the distance function $d(x,y)$ on $M$ (with respect to the smooth background metric $g$). When $y$ is sufficiently close to $x$, and if $\gamma$ is a unit-speed minimizing geodesic connecting $x$ to $y$, then
\[
\nabla_x d(x,y) = -\gamma'(0) \text{ and } \nabla_y d(x,y) = \gamma'(d(x,y)).
\]
Thus, if $P(x,y): T_y M \rightarrow T_x M$ denotes the parallel transport map along the unique shortest geodesic connecting $x$ and $y$, then 
\[
\nabla_x d(x,y) = -P(x,y)(\nabla_y d(x,y)).
\]
$P(x,y)$ can also be considered as an element of $T_{(x,y)} M \times M$ by letting $P(x,y)(X, Y) = g(X,P(x,y)Y)$. $P$ is then a covector field on an open neighborhood of the diagonal of $M \times M$. 

In the sequel, we represent $P$ in local coordinates by using two indices (which can be casually raised or lowered using the metric $g$): the first index refers to the $x$-factor and the second stands for the $y$-factor. For example, as a transformation of $T_y M$ into $T_x M$, we have
\[
P(x,y) = P^i{}_j(x,y)\, \partial_{x^i} \otimes dy^j,
\]
while, as a covector field, we have
\[
P(x,y) = P_{ij}(x,y)\,dx^i\,dy^j.
\]

We make a few observations:
\begin{enumerate}[(P1)]
\item $P(x,x) = Id$.

\item For any compact subset $K$ of $M$, there exists $\delta = \delta(K)$ such that $P$ is smooth in $\{(x,y) \in K \times K: d(x,y) < \delta\}$. 

\item $\nabla_x P(x,x) = 0$ and $\nabla_y P(x,x) = 0$. To see this, pick any geodesic $\gamma(t)$ emanating from $x$ (so that $\gamma(0) = x$). Then $P(x,\gamma(t))$ is parallel along $\gamma$, i.e.  $\nabla_{\gamma'(t)} P(x,\gamma(t)) = 0$. As $\gamma'(0)$ was chosen arbitrarily, this gives $\nabla_y P(x,x) = 0$. Likewise, $\nabla_x P(x,x) = 0$.

\item It holds that
\begin{equation}
g^{kl}(y) = g^{ij}(x)\,P_j{}^k(x,y)\,P_i{}^l(x,y).
	\label{Eq:QOrt}
\end{equation}
To see this, take any covector $Y \in T^*_y M$ and let $X = (P(x,y)Y^\sharp)_\flat \in T^*_x M$. Then $X_i = P_i{}^l(x,y)\,Y_l$ and so
\[
g^{kl}(y)\,Y_k\,Y_l = g^{ij}(x)\,X_i\,X_j = g^{ij}(x)\,P_j{}^k(x,y)\,P_i{}^l(x,y)\,Y_k\,Y_l.
\]
Since $Y$ is arbitrary, this implies the asserted identity.

Note that \eqref{Eq:QOrt} implies that
\[
[g_{rk}(y)\, g^{ij}(x)\,P_j{}^k(x,y)]\,P_i{}^l(x,y) = \delta_r^l,
\]
which further implies that
\[
P_t{}^r(x,y) [g_{rk}(y)\, g^{ij}(x)\,P_j{}^k(x,y)] = \delta_t^k
\]
and
\begin{equation}
g_{tj}(x) = g_{rk}(y)\,P_t{}^r(x,y)\,P_j{}^k(x,y).
	\label{Eq:QOrtX}
\end{equation}
\end{enumerate}

\begin{proof}[Proof of Proposition \ref{Prop:LipWeakStab}.] {\it Step 1:} We start with a decomposition of the leading order term in $\Ric(e^{2\bar f_\vareps}g)$. We compute
\begin{align*}
\nabla_{x^i} \nabla_{x^j} \bar f_\vareps(x) 
	&= \int_M \nabla_{x^i} \nabla_{x^j} \varrho_\vareps(d(x,y))\,f(y)\,dv_g(y)\\
	&= -\int_M  \nabla_{x^i}[P_j{}^k(x,y) \nabla_{y^k} \varrho_\vareps(d(x,y))]\,f(y)\,dv_g(y)\\
	&= -\int_M  P_j{}^k(x,y) \nabla_{x^i} \nabla_{y^k} \varrho_\vareps(d(x,y))\,f(y)\,dv_g(y)\\
		&\qquad\qquad - \int_M \nabla_{x^i}P_j{}^k(x,y) \,\nabla_{y^k} \varrho_\vareps(d(x,y))\,f(y)\,dv_g(y)\\
	&= \int_M  P_j{}^k(x,y)  \nabla_{x^i}  \varrho_\vareps(d(x,y))\,\nabla_{y^k} f(y)\,dv_g(y)\\
		&\qquad\qquad + \int_M \nabla_{y^k} P_j{}^k(x,y) \nabla_{x^i}  \varrho_\vareps(d(x,y))\, f(y)\,dv_g(y)\\
		&\qquad\qquad - \int_M  \nabla_{x^i}P_j{}^k(x,y) \,\nabla_{y^k} \varrho_\vareps(d(x,y))\,f(y)\,dv_g(y)\\
	&=: T^{(1)}_{ij}(x) + T^{(2)}_{ij}(x) + T^{(3)}_{ij}(x).
\end{align*}
An analogous calculation also gives
\begin{align*}
\nabla_{x^i} \nabla_{x^j} Z_\vareps(x) 
	&= \int_M \nabla_{y^k} P_j{}^k(x,y) \nabla_{x^i}  \varrho_\vareps(d(x,y))\,dv_g(y)\\
		&\qquad\qquad - \int_M  \nabla_{x^i}P_j{}^k(x,y) \,\nabla_{y^k} \varrho_\vareps(d(x,y))\,dv_g(y),
\end{align*}
where $Z_\vareps$ is as defined in \eqref{Eq:05I18-ZepsDef}. Keeping in mind that $\nabla_x P(x,x) = 0$ and $\nabla_y P(x,y) = 0$, we thus deduce that 
\begin{align*}
|T^{(2)}_{ij}(x) + T^{(3)}_{ij}(x)|_g
 	&= \Big|\int_M \nabla_{y^k} P_j{}^k(x,y) \nabla_{x^i}  \varrho_\vareps(d(x,y))\, [f(y) - f(x)]\,dv_g(y)\\
		&\qquad\qquad - \int_M  \nabla_{x^i}P_j{}^k(x,y) \,\nabla_{y^k} \varrho_\vareps(d(x,y))\,[f(y) - f(x)]\,dv_g(y)\\
		&\qquad\qquad + f(x)\nabla_{x^i} \nabla_{x^j} Z_\vareps(x)  \Big|_g\\
	&\leq o(1) \|f\|_{C^{0,1}(\Omega)},
\end{align*}
where, here and below, $o(1)$ denotes some constant such that $\lim_{\vareps \rightarrow 0} o(1) = 0$.

We also have
\begin{align*}
T^{(1)}_{ij}(x)
	&= \int_M  P_j{}^k(x,y)  \nabla_{x^i}  \varrho_\vareps(d(x,y))\,\nabla_{y^k} f(y)\,dv_g(y)\\
	&= -\int_M  P_j{}^k(x,y)  P_i^l(x,y)\, \nabla_{y^l}  \varrho_\vareps(d(x,y))\,\nabla_{y^k} f(y)\,dv_g(y)\\
	&= -\int_M  \nabla_{y^l}[P_j{}^k(x,y)  P_i^l(x,y)\,  \varrho_\vareps(d(x,y))]\,\nabla_{y^k} f(y)\,dv_g(y)\\
		&\qquad\qquad + T^{(4)}_{ij}(x)\\
	&=: T^{(0)}_{ij}(x) + T^{(4)}_{ij}(x)		
		,
\end{align*}
where $|T^{(4)}_{ij}(x)|_g \leq o(1) \|f\|_{C^{0,1}(\Omega)}$.

We thus have
\[
-(n-2)\nabla_g^2\bar f_\vareps - \Delta_g \bar f_\vareps\,g 
	\geq -(n-2) T^{(0)} - \tr_g(T^{(1)})\,g - o(1)\,\|f\|_{C^{0,1}(\Omega)}.
\]

Since $\nabla \bar f_\vareps \rightarrow \nabla f$ in $L^p_{loc}(\Omega)$ (and uniformly if $f \in C^1(\bar\Omega)$), to establish the result, it suffices to show that
\begin{equation}
-(n-2) T^{(0)} - \tr_g(T^{(1)})\,g \geq \tilde F
	\label{Eq:05I18-QRicespLB}
\end{equation}
where the $(0,2)$-tensor $\tilde F$ is defined by
\begin{align*}
\tilde F_{ij}(x)
	&= \int_M F_{kl}(y)\,P_j{}^k(x,y)  P_i^l(x,y)\,  \varrho_\vareps(d(x,y))\,dy,\\
F_{ij} 
	&= -(n-2)\nabla_i f \, \nabla_j f + (n-2)|df|_g^2\,g_{ij} - \Ric_{ij}(g) + (n-1)k\,e^{2f}\,g_{ij}.
\end{align*}

Let $a$ be some non-negative symmetric $(2,0)$-tensor $a$ with compact support in $\Omega$. Define a $(2,0)$-tensor $a_\vareps$ defined by
\[
a_{\vareps}^{kl}(y) = \int_M  a^{ij}(x)  P_j{}^k(x,y)  P_i^l(x,y)\,  \varrho_\vareps(d(x,y))\,dv_g(x).
\]
Then $a_\vareps$ is symmetric and non-negative, as it holds for any covector $V \in T_y^*M$ that
\begin{align*}
a_\vareps(y)(V,V) 
	&= \int_M  a^{ij}(x)  P_j{}^k(x,y) V_k  P_i^l(x,y)\,V_l  \varrho_\vareps(d(x,y))\,dv_g(x)\\
	&= \int_M  \underbrace{a(x)(  P(x,y)(V), P(x,y)(V))}_{\geq 0} \varrho_\vareps(d(x,y))\,dv_g(x)
		\geq 0.
\end{align*}

We have
\begin{align*}
&\int_M a^{ij}(x)\,T^{(0)}_{ij}(x)\,dv_g(x)\\
	&\qquad= -\int_M a^{ij}(x) \int_M  \nabla_{y^l}[P_j{}^k(x,y)  P_i^l(x,y)\,  \varrho_\vareps(d(x,y))]\,\nabla_{y^k} f(y)\,dv_g(y) \,dv_g(x)\\
	&\qquad= -\int_M \nabla_{y^k} f(y) \int_M  \nabla_{y^l}[a^{ij}(x)  P_j{}^k(x,y)  P_i^l(x,y)\,  \varrho_\vareps(d(x,y))]\,\,dv_g(x) \,dv_g(y)\\
	&\qquad = -\int_M \nabla_{y^k} f(y) \nabla_{y^l} a_{\vareps}^{kl}(y) \,dv_g(y),
\end{align*}
In addition, since 
\begin{align*}
\tr_g(a_\vareps)(y) 
	&= g_{kl}(y) a_{\vareps}^{kl}(y) = \int_M  a^{ij}(x)  \underbrace{g_{kl}(y) P_j{}^k(x,y)  P_i^l(x,y)}_{\stackrel{\eqref{Eq:QOrtX}}{=} g_{ij}(x)}\,  \varrho_\vareps(d(x,y))\,dv_g(x)\\
	&= \int_M  \tr_g(a)(x)  \varrho_\vareps(d(x,y))\,dv_g(x),
\end{align*}
we also have that
\begin{align*}
&\int_M \tr_g(a)(x) \tr_g(T^{(1)})(x)\,dv_g(x)\\
	&\quad= -\int_M \tr_g(a)(x)\,\int_M  \underbrace{g^{ij}(x) P_j{}^k(x,y)  P_i^l(x,y)}_{\stackrel{\eqref{Eq:QOrt}}{=}g^{kl}(y)}\,  \nabla_{y^l}\varrho_\vareps(d(x,y))\,\nabla_{y^k} f(y)\,dv_g(y) \,dv_g(x)\\
	&\quad= -\int_M \tr_g(a)(x)\,\int_M  g^{kl}(y)  \nabla_{y^l}\varrho_\vareps(d(x,y))\,\nabla_{y^k} f(y)\,dv_g(y) \,dv_g(x)\\
	&\quad= -\int_M g^{kl}(y) \nabla_{y^k} f(y) \nabla_{y^l}\Big\{ \int_M \tr_g(a)(x)    \varrho_\vareps(d(x,y))\,dv_g(x)\Big\} \,dv_g(y)\\
	&\quad = -\int_M \nabla_{y^k} f(y) \nabla^{y^k}\tr_g(a_\vareps)(y)  \,dv_g(y).
\end{align*}
It hence follows that
\begin{align*}
&\int_M a^{ij}(x)\,[-(n-2) T^{(0)}_{ij}(x) -  \tr_g(T^{(1)})\,g_{ij}] \,dv_g(x)\\
	&\qquad = \int_M \nabla_{y^k} f(y) \Big[(n-2) \nabla_{y^l} a_{\vareps}^{kl}(y) +  \nabla^{y^k}\tr_g(a_\vareps)(y)\Big] \,dv_g(y)
\end{align*}
Now since $\Ric(e^{2f}g) \geq (n-1)k$ in the weak sense (by Proposition \ref{Prop:Vis=>Weak}) and $a_\vareps$ is non-negative definite, we arrive at
\[
\int_M a^{ij}(x)\,[-(n-2) T^{(0)}_{ij}(x) -  \tr_g(T^{(1)})\,g_{ij}] \,dv_g(x)
	\geq \int_M F_{ij}(y)\,a_\vareps^{ij}(y)\,dv_g(y),
\]
from which \eqref{Eq:05I18-QRicespLB} is readily seen. This completes the proof.
\end{proof}

\subsection{Volume comparison}

We are now ready to give the proof of the relative volume comparison theorem for continuously conformally metrics with lower Ricci bounds.

\begin{proof}[Proof of Theorem \ref{Thm:RelVolComp}] By Propositions \ref{Prop:InfConvStab} and \ref{Prop:LipWeakStab}, there exists a sequence of smooth functions $\{\bar f_\vareps\}$ such that, as $\vareps \rightarrow 0$, $\bar f_\vareps \rightarrow f$ locally uniformly in $\Omega$ and $\{\bar f_\vareps\}$ satisfies an integral Ricci lower bound
\[
\lim_{\vareps \rightarrow 0} \int_{\omega} \Big\{ \max\big(-\lambda_1(\Ric(e^{2\bar f_\vareps}g)) + (n-1)k, 0\big)\Big\}^p\,dv_g = 0
\]
for any open $\omega \Subset \Omega$ and any $1 \leq p < \infty$.

Let
\[
\Lambda_\vareps(\omega,p) = \int_{\omega} \Big\{ \max\big(-\lambda_1(\Ric(e^{2\bar f_\vareps}g)) + (n-1)k, 0\big)\Big\}^p\,dv_g.
\]
Then, for $p > \frac{n}{2}$, the relative volume comparison theorem of Petersen and Wei \cite[Theorem 1.1]{PetersenWei-GAFA97} (see also \cite{Wei-SDG07}) implies for $0 < r < R$ that
\[
\Big(\frac{Vol_{e^{2\bar f_\vareps}g}(B_{e^{2\bar f_\vareps}g}(p,R))}{v(n,k,R)}\Big)^{\frac{1}{2p}} - \Big(\frac{Vol_{e^{2\bar f_\vareps}g}(B_{e^{2\bar f_\vareps}g}(p,r))}{v(n,k,r)}\Big)^{\frac{1}{2p}} \leq C(R)\, \Lambda_\vareps(\omega,p)^{\frac{1}{2p}}.
\]
(Here we assume $R \leq \frac{\pi}{2\sqrt{k}}$ if $k > 0$.) Sending $\vareps \rightarrow 0$ we obtain the first conclusion. 

We turn to the second conclusion. By another theorem of Petersen and Wei \cite[Theorem 1.5]{PetersenWei-TrAMS01}, there is a map $\phi: B_{e^{2f}g}(p,r) \rightarrow \SSphere^n_k$ which preserves the distance function. We need to show that $\phi$ and $f$ are smooth. 

We represent $\phi(B_{e^{2f}g}(p,r))$ as a ball $B(0,\tilde r) \subset \RR^n$ equipped with a conformally flat metric $g_{can} = e^{2F}g_{flat}$ where $g_{flat}$ is the flat metric on $\RR^n$ and $F$ is a smooth function. Let $\{x^1, \ldots, x^n\}$ be a local coordinate system on $M$ relative to which $g$ is smooth. Let $\{y^1, \ldots, y^n\}$ denote a standard coordinate system on $\RR^n$.

Observe that $\phi$ considered as a map from $(B_{e^{2f}g}(p,r),g)$ into $(B(0,\tilde r), g_{flat})$ is locally Lipschitz continuous (since $f$ is locally bounded). Hence $\phi$ is differentiable almost everywhere. Likewise, $\psi := \phi^{-1}$ is differentiable almost everywhere. 

We claim that $e^{2F}g_{flat} = \psi^*(e^{2f}g)$, i.e.
\begin{equation}
e^{2F(y)} \delta_{ij} = e^{2f(\psi(y))} g_{kl}(\psi(y)) \frac{\partial\psi^k}{\partial y^i}(y) \frac{\partial\psi^l}{\partial y^j}  (y) \text{ a.e. in } B(0,\tilde r).
	\label{Eq:11II18-E1}
\end{equation}

We will use the following formula (see e.g. \cite[Theorem 2.7.6]{Burago2Ivanov}) for the length of a Lipchitz curve $\gamma: [a,b] \rightarrow X$ in a metric space $(X,d)$ where the distance function $d$ is generated by a metric $e^{2u} h$ where $u$ is continuous and $h$ is smooth:
\[
\textrm{Length}_d(\gamma([a,b])) = \int_a^b e^{u(\gamma(t))} |\gamma'(t)|_h\,dt.
\]
(Here we are using that 
\[
 \lim_{\delta \rightarrow 0} \frac{d(\gamma(t + \delta),\gamma(t))}{\delta} =  e^{u(\gamma(t))} \lim_{\delta \rightarrow 0} \frac{d_h (\gamma(t + \delta),\gamma(t))}{\delta} = e^{u(\gamma(t))} |\gamma'(t)|_h
 \]
 at points where $\gamma$ is differentiable.)

We note that, since $\psi$ preserves the distance, it preserves lengths of curves. Hence if $\gamma: [a,b] \rightarrow B(0,\tilde r)$ is a Lipschitz curve, then
\begin{align*}
\int_a^b e^{F \circ \gamma}|\gamma'(t)|_{g_{flat}}\,dt 
	&= \textrm{Length}_{e^{2F}g_{flat}}(\gamma([a,b])) \\
	&= \textrm{Length}_{e^{2f}g} (\psi \circ \gamma([a,b]))
	= \int_a^b  e^{f \circ \psi \circ \gamma}(t) \Big| \frac{d}{dt} (\psi \circ \gamma) \Big|_g(t)\,dt\\
	&= \int_a^b e^{f \circ \psi \circ \gamma}(t) \Big(g^{kl} \circ \psi \circ \gamma \frac{d}{dt} (\psi^k \circ \gamma)\, \frac{d}{dt} (\psi^l \circ \gamma)\Big)^{1/2}\,dt.
\end{align*}

Now, for each $i \in \{1, \ldots, n\}$, consider the family of curves
\[
\gamma_{y_1, \ldots, \hat y^i, \ldots, y_n}(t) = (y_1, \ldots, t y^i, \ldots, y_n),
\]
where the hat above $y^i$ indicates that this entry is absent. We then have
\begin{align*}
\int_a^b e^{F \circ \gamma_{y_1, \ldots, \hat y^i, \ldots, y_n}}\,dt 
	&= \int_a^b e^{f \circ \psi} \Big(g_{kl}(\psi) \frac{\partial\psi^k}{\partial y^i} \frac{\partial\psi^l}{\partial y^i}\Big)^{1/2}\Big|_{y = \gamma_{y_1, \ldots, \hat y^i, \ldots, y_n}} \,dt
\end{align*}
for almost all $(y_1, \ldots, \hat y^i, \ldots, y_n) \in \RR^{n-1}$ and for all $a, b$ such that $\gamma_{y_1, \ldots, \hat y^i, \ldots, y_n}([a,b]) \subset B(0,\tilde r)$. This implies that, for every $i \in \{1, \ldots, n\}$, 
\[
e^{2F(y)} = e^{2f(\psi(y))} g_{kl}(\psi(y)) \frac{\partial\psi^k}{\partial y^i}(y) \frac{\partial\psi^l}{\partial y^i}  (y) \text{ a.e. in } B(0,\tilde r). 
\]
Similarly, by considering family of curves tangential to $\partial_{y^i} + \partial_{y^j}$, we have, for every $i, j \in \{1, \ldots, n\}$,
\begin{multline*}
2e^{2U(y)} = e^{2f(\psi(y))} g_{kl}(\psi(y)) \Big(\frac{\partial\psi^k}{\partial y^i}(y) + \frac{\partial\psi^k}{\partial y^j}  (y)\Big) \times\\
	\times \Big(\frac{\partial\psi^l}{\partial y^i}(y) + \frac{\partial\psi^l}{\partial y^j}  (y)\Big) \text{ a.e. in } B(0,\tilde r). 
\end{multline*}
The claim \eqref{Eq:11II18-E1} follows from the above two equations.

For $D \subset \RR^n$ and $u \in W^{1,n}(D)$, consider the functional 
\[
I[u;D] = \int_{D} |\nabla_{flat} u|_{g_{flat}}^n\,dv_{flat}  = \int_{D} |\nabla_{g_{can}} u|_{g_{can}}^n\,dv_{g_{can}}.
\]
Similarly, for $D \subset \Omega$ and $u \in W^{1,n}(D,g)$, consider
\[
J[u; D] = \int_D |\nabla_{e^{2f} g} u|_{e^{2f} g}^n\,dv_{e^{2f} g}  = \int_{D} |\nabla_{g} u|_{g}^n\,dv_{g}.
\]

Observe that, by convexity, for each $1 \leq i \leq n$, the function $y^i$ on $B(0,\tilde r) \subset \RR^n$ satisfies, for $D \subset B(0,\tilde r)$, that 
\[
I[y^i;D] \leq I[u;D] \text{ for all } u \in W^{1,n}(D) \text{ such that }u = y^i \text{ on } \partial D.
\]
Noting that $y^i(y) = \phi^i(\psi(y))$, and using the fact that the change of variable formula holds for Lipschitz transformation (see e.g. \cite[p. 99]{EvansGariepy92}), we have, for $D \subset B_{e^{2f}g}(p,r)$, 
\[
J[\phi^i; D] \leq J[u; D] \text{ for all } u \in W^{1,n}(D,g) \text{ such that }u = \phi^i \text{ on } \partial D.
\]
It follows that $\phi^i$ satisfies
\begin{equation}
\textrm{div}_g(|\nabla_g \phi^i|_g^{n-2}\nabla_g \phi^i) = 0 \text{ in } B_{e^{2f}g}(p,r).
	\label{Eq:7II18-E1}
\end{equation}
Noting also that $|\nabla y^i|_{g_{flat}} = 1$, we can find $C > 1$ such that  $C^{-1} < |\nabla_g \phi^i|_g < C$ in $B_{e^{2f}g}(p,r)$. It follows that equation \eqref{Eq:7II18-E1} is a uniformly elliptic quasilinear equation. A regularity result of Ladyzhenskaya and Uraltseva for quasilinear and uniformly elliptic (scalar) equations in divergence form (\cite[Chapter 4]{LadyzhenskayaUraltseva68}) implies that $\phi^i$ belongs to $W^{2,2}_{loc}$ and $C^{1,\alpha}_{loc}$ for some $\alpha \in (0,1)$. (The $C^{1,\alpha}$ regularity also follows from \cite{Uraltseva68-LOMI, Uhlenbeck77-ActaM, Evans82-JDE, Lewis80-PrAMS, DiBenedetto83-NA} where $|\nabla_g \phi^i|_g$ is allowed to vanish.) We then recast equation \eqref{Eq:7II18-E1} in non-divergence form
\[
A^{kl}(\nabla_g \phi)\,\nabla_k \nabla_l \phi^i = 0,
\]
which is understood in the almost everywhere sense and where the coefficients $A^{kl}$ is uniformly elliptic. Now, as a function of $x$, $A^{kl}(\nabla_g \phi(x))$ is $C^{\alpha}$ continuous, and so elliptic regularity implies that $\phi^i$ is $C^{2,\alpha}$. The smoothness of $\phi^i$ follows from bootstrapping.  Recalling that $e^{2f}g = \phi^*(g_{can})$, we deduce that $f$ is smooth and conclude the proof.
\end{proof}

\section{Nonexistence of Green's functions for $\mu_\Gamma^+ \leq 1$}\label{Sec:NoExistence}

In this section, we prove part (ii) of Theorem \ref{thm:MainThm}. In fact, we have:

\begin{theorem}\label{thm:M(ii)}
Let $(M,g)$ be an $n$-dimensional smooth compact Riemannian manifold with $n \geq 3$. Assume that $\Gamma$ satisfies \eqref{G1}, \eqref{G2} and that $\lambda(A_g) \in \Gamma$ in $M$. Let $S = \{p_1, \ldots, p_m\}$ be a non-empty finite subset of distinct points of $M$ and $c_1, \ldots, c_m \in (0,\infty)$. If $\mu_\Gamma^+ \leq 1$, then the following are equivalent
\begin{enumerate}[(i)]
\item there exists a function $u \in C^0(M\setminus S)$ such that 
\begin{align}
&\lambda(A_{g_u}) \in \bar\Gamma \text{ and }  u > 0 \text{ in } M \setminus \{p_1, \ldots, p_m\},\label{Eq:DegEq1Super}\\
&\lim_{x \rightarrow p_i} d_g(x,p_i)^{n-2} u(x) = c_i,\qquad i = 1, \ldots, m, \label{Eq:DegAsSuper}
\end{align}
where $d_g$ denotes the distance function with respect to the metric $g$;
\item $(M,g)$ is conformal to the standard sphere and $m = 1$. 
\end{enumerate}
\end{theorem}

\begin{proof}
It is clear that (ii) implies (i). Conversely, assume that (i) holds. Then, by the relative volume comparison theorem (Theorem \ref{Thm:RelVolComp}), for any $p \in M \setminus S$, the function
\[
r \mapsto \frac{Vol_{g_u}(B_{g_u}(p,r))}{\omega_n\,r^n}
\]
is non-increasing, where $\omega_n$ is the volume of the unit $n$-dimensional Euclidean unit ball. On the other hand, as $r \rightarrow 0$, the above function tends to $1$, and, as $r \rightarrow \infty$, it tends to $m$ (thanks to \eqref{Eq:DegAsSuper}). It follows that $m = 1$ and that $Vol_{g_u}(B_{g_u}(p,r)) = \omega_n\,r^n$ for all $r > 0$. By the rigidity part of the relative volume comparison theorem, we have that $u$ is smooth and $(M\setminus S, g_u)$ is isometric to the Euclidean space $\RR^n$. We then proceed as in \cite[Section 7.6]{GV07}: The metric $g$ is conformally flat on $M \setminus S$ and so is locally conformally flat on $M$ by the vanishing of the Weyl tensor for $n \geq 4$ and of the Cotton tensor for $n = 3$. In addition, $M$, being a one-point compactification of $M \setminus S$, is homeomorphic to $\SSphere^n$, and hence is simply connected. A theorem of Kuiper \cite[Theorem 6]{Kuiper49} then implies that $(M, g)$ is conformally equivalent to the standard sphere. 
\end{proof}


\section{Existence and uniqueness of Green's functions for $\mu_\Gamma^+ > 1$}\label{Sec:ExistenceUniq}

In this section we prove part (i) of Theorem \ref{thm:MainThm}. For simplicity, we will only present the proof in the case where $S$ consists of a single point and $c_1 = 1$. The proof can be easily adapted to treat the general case.

\subsection{Non-degenerate elliptic Dirichlet boundary value problems}\label{Sec:DBVP}

Let $\Gamma$ satisfy \eqref{G1}, \eqref{G2} and $f$ satisfy \eqref{f1}-\eqref{f4}. It is easily seen that equation \eqref{Eq:DegEq1} is the same as 
\[
f\big(\lambda(A_{g_u})\big) = 0 \text{ on } M \setminus \{p_1, \ldots, p_m\}.
\]
We will eventually regularize this equation by replacing the right hand side by small positive constants.

\begin{theorem}\label{Thm:GuanGen}
Let $n \geq 3$ be an integer and $(\bar N,g)$ be an $n$-dimensional smooth compact Riemannian manifold with non-empty smooth boundary $\partial N$. Assume that $(f,\Gamma)$ satisfies \eqref{G1}-\eqref{G2}, \eqref{f1}-\eqref{f4}. Let $\psi \in C^\infty(\bar N \times \RR)$, $\psi > 0$ and $\varphi \in C^\infty(\partial N)$. Assume that there exists a function $\bar u \in C^\infty(\bar N)$ such that $\bar u \equiv \varphi$ on $\partial N$ and
\[
f(\lambda(A_{g_{\bar u}})) \geq \psi(\cdot, \bar u) \text{ in } \bar N.
\]
Then, there exists a solution $u \in C^\infty(N) \cap C^{0,1}(\bar N)$ (with $u \leq \bar u$ in $\bar N$) to the boundary value problem
\begin{align}
&f\big(\lambda(A_{g_u})\big) = \psi(\cdot, u) \text{ in }  N,\label{Eq:NondegEq}\\
&u = \varphi \text{ on } \partial N.\label{Eq:NdBC}
\end{align}
Moreover, there exists a constant $C > 0$ depending only on $(\bar N,g)$, $(f,\Gamma)$, $\psi$, $\varphi$, $\|\ln \bar u\|_{C^3(\bar N)}$ and $\lambda(A_{g_{\bar u}})$ such that 
\[
\|\ln u\|_{C^{0,1}(\bar N)} \leq C,
\]
and for every compact subset $K$ of $N$ and every $l \geq 2$, there exists $C_{K,l}$ depending only on $K$, $l$, $(\bar N,g)$, $(f,\Gamma)$,  $\psi$, $\varphi$, $\|\ln \bar u\|_{C^3(\bar N)}$ and $\lambda(A_{g_{\bar u}})$ so that
\[
\|\ln u\|_{C^l(K)} \leq C_{K,l}.
\]
\end{theorem}

When $(f,\Gamma) = (\sigma_k^{\frac{1}{k}},\Gamma_k)$, the result was proved in Guan \cite{Guan07-AJM}. In fact, in this case, the proof therein yields $C^2$-estimate up to the boundary. We chose to forgo such estimate in full generality as it is not needed for our current purpose. We instead circumvent the issue by ``opening up'' $\Gamma$ to larger cones $\Gamma_t$ where a double normal derivative estimate for $\Gamma_t$ can be obtained fairly easily. The procedure in \cite{Guan07-AJM} can then be applied to prove the existence of solutions corresponding to those cones $\Gamma_t$. Letting $\Gamma_t$ converge back to $\Gamma$,  we obtain Theorem \ref{Thm:GuanGen} above.

\begin{proof}
Replacing $g$ by $g_{\bar u}$ if necessary, we may assume that $\lambda(A_g) \in \Gamma$. Let $\underline{u}$ be the solution to 
\begin{align*}
L_g \underline{u} & = 0 \text{ in } N,\\
\underline{u}	& = \varphi \text{ on } \partial N.
\end{align*}
By \eqref{G2}, $\underbar{u}$ is a subsolution to \eqref{Eq:NondegEq}. In particular $\underline{u} \leq \bar u$. We will construct a solution to \eqref{Eq:NondegEq}-\eqref{Eq:NdBC} which satisfies $\underline{u} \leq u \leq \bar u$. We will argue according to whether $(1,0,\ldots, 0) \in \Gamma$ or not.

\medskip
\noindent{\it Step 1:} Assume for the moment that $(1,0,\ldots, 0) \in \Gamma$. We adapt the argument in \cite{Guan07-AJM} to the case at hand.

By mean of a degree theoretic argument (and Evans-Krylov estimates), it suffices to show that, there exists a constant $C$ such that if $u$ is a solution to \eqref{Eq:NondegEq}-\eqref{Eq:NdBC} satisfying $\underline{u} \leq u \leq \bar u$ then
\begin{equation}
\|u\|_{C^2(N)} \leq C.
	\label{Eq:04VI19-E1}
\end{equation}

Since $(0, \ldots, 0, 1) \in \Gamma$ and $f$ is homogeneous of degree one, there exists $\delta = \delta(f,\Gamma) > 0$ such that for every compact set $E \subset \Gamma$, there exists $\bar R = R(\delta, E) > 0$ such that, for all $\lambda = (\lambda_1, \ldots, \lambda_n) \in E$ and $R > \bar R$,
\begin{equation}
f(\lambda_1, \ldots, \lambda_{n-1},\lambda_n + R) = R f(\frac{1}{R}\lambda + (0,\ldots, 0, 1)) > R\delta > 0.
	\label{Eq:BG1.13}
\end{equation}
(This implies \cite[eq. (1.13)]{Guan07-AJM}.) Also, we claim that 
\begin{equation}
\sum_{i=1}^n f_{\lambda_i}(\lambda) \geq f\big(1, \ldots, 1\big) > 0 \text{ in } \Gamma.
	\label{Eq:29-1}
\end{equation}
(This is equivalent to \cite[eq. (1.10)]{Guan07-AJM}.) To see this, let $e = (1, \ldots, 1)$. For every $\lambda \in \Gamma$ and $\mu > 0$, we have $\lambda + \mu e \in \Gamma$ due to \eqref{G1}-\eqref{G2}. The concavity of $f$ then gives $f(\lambda + \mu e) \leq f(\lambda) + \mu \sum_i f_{\lambda_i}(\lambda)$. Dividing by $\mu$ and letting $\mu \rightarrow \infty$, we obtain \eqref{Eq:29-1} in view of the homogeneity of $f$.

In view of \eqref{Eq:29-1}, the proof of \cite[Theorem 3.3 and Theorem 3.4]{Guan07-AJM} can be applied directly to the present setting yielding
\[
\max_N |\nabla \ln u| \leq C_1 \text{ and } \max_N |\nabla^2 \ln u| \leq C_2,
\]
where $C_1$ depends on $(M,g)$, $(f,\Gamma)$, $\max_N |\ln u|$, $\max_{\partial N} |\nabla \ln u|$ and  $C_2$ depends on $(M,g)$, $(f,\Gamma)$, $\max_N |\ln u|$, $C_1$ and $\max_{\partial N} |\nabla^2 \ln u|$. (To dispel confusion, note that the function $u$ appearing in \cite{Guan07-AJM} is $\ln\frac{1}{u}$ in our present setting. Also, the parameters $s$ and $t$ therein are taken to be $1$.) As $u$ is pinched between $\underline{u}$ and $\bar u$, $\max_{\partial N} |\nabla \ln u|$ is bounded in terms of $|\nabla \varphi|$, $|\partial_\nu \bar u|$ and $|\partial_\nu \underline{u}|$, where $\nu$ is the unit normal to $\partial N$. Thus, to establish \eqref{Eq:04VI19-E1}, it suffices to show that
\begin{equation}
|\nabla^2 u| \leq C \text{ on } \partial N,
	\label{Eq:Eq2.1}
\end{equation}
where $C$ depends on $(\bar N,g)$, $(f,\Gamma)$, $\psi$, $\varphi$, $\|u\|_{C^1(\bar N)}$ and $\lambda(A_{g_{\bar u}})$.

For $x_0 \in \partial N$, let $e_1, \ldots, e_n$ be an orthonormal frame about $x_0$ obtained by parallel transporting an orthonormal local frame $e_1, \ldots, e_{n-1}$ on $\partial N$ and the inward pointing unit normal $e_n$ to $\partial N$ along geodesics perpendicular to $\partial N$.

Let $\mathcal{L} = \sum_{ij} F^{ij} \nabla_i \nabla_j$ be the principal part of the linearized operator for $f(\lambda(A_{g_u}))$ at $u$. Using \eqref{Eq:BG1.13}, one can check that \cite[Lemma 2.2]{Guan07-AJM} holds in the present setting: For any $B > 0$, there exist small positive constants $\mu$ and $\delta$ and a large positive constant $N$ such that the function
\[
v = -\ln \frac{u}{\bar u} - \frac{1}{2}\Big(\ln \frac{u}{\bar u}\Big)^2 + \mu d( \cdot, \partial N) - \frac{1}{2} N d^2( \cdot, \partial N)
\]
satisfies 
\[
\mathcal{L}[v] \leq - B - \beta \sum_i F^{ii}.
\]
We can now follow the proof of \cite[eqs. (2.10), (2.12)]{Guan07-AJM} to obtain
\[
|\nabla_{ij} u(x_0)| + |(A_{g_u})_{ij}(x_0)| \leq C_0 \text{ provided } (i,j) \neq (n,n).
\]
Since $u$ is super-harmonic (with respect to the conformal Laplacian), this implies that
\[
\nabla_{nn} u(x_0) \geq - C \text{ and } (A_{g_u})_{nn}(x_0) \geq -C.
\]

It remains to give an upper bound for $\nabla_{nn} u(x_0)$, where our argument differs from (and is much easier than) that in \cite{Guan07-AJM} (where some algebraic properties of the $\sigma_k$-equation play more of a role).  Since $(1, 0, \ldots, 0) \in \Gamma$, there exists some $C_1 > 0$ such that if $|a_{ij}| < \frac{1}{C_1}$ for $(i,j) \neq (n,n)$ and $a_{nn} = 1$, then $\lambda((a_{ij})) \in \Gamma$ and $f(\lambda((a_{ij}))) \geq \frac{1}{C_1}$.

If $(A_{g_u})_{nn}(x_0) \leq C_0\,C_1$, we are done. Otherwise, we have
\begin{align*}
\psi(x,u) 
	&= f(\lambda(A_{g_u}(x_0)))
		 = (A_{g_u})_{nn}(x_0)f\Big(\lambda\Big(\frac{1}{(A_{g_u})_{nn} (x_0)} A_{g_u}(x_0)\Big)\Big)\\
	& \geq (A_{g_u})_{nn}(x_0)\,\frac{1}{C_1}.
\end{align*}
This implies that $(A_{g_u})_{nn}(x_0) \leq C_1\,\psi(x_0,u(x_0)) \leq CC_1$. We have thus established \eqref{Eq:Eq2.1}, and thus established the theorem when $(1,0, \ldots, 0) \in \Gamma$.

\medskip
\noindent{\it Step 2:} We now return to the general case where $(1,0,\ldots, 0)$ may or may not belong to $\Gamma$. For $t \in [\frac{1}{2},1]$, define 
\begin{align*}
\Gamma_t
	&:=\{\lambda\in \Bbb R^n\ |\
t\lambda+(1-t)\sigma_1(\lambda)e\in \Gamma\}, \quad
\mbox{where}\ e=(1,\cdots, 1),\\
f_t(\lambda)
	&=f(t\lambda+(1-t)\sigma_1(\lambda)e).
\end{align*}
It was proved in \cite{LiLi03} that $(f_t, \Gamma_t)$ 
also satisfies (\ref{f1})-(\ref{f3}). 

Note that $(1,0,\ldots, 0) \in \Gamma_t$ for $t < 1$ as $(1, 1-t, \ldots, 1 -t) \in \Gamma_n \subset \Gamma$. Furthermore, we have $f_t(\lambda) \geq f(t\lambda) = tf(\lambda)$ and so $\bar u$ satisfies
\[
f_t(\lambda(A_{g_{\bar u}})) \geq t\psi(\cdot, \bar u) \text{ in } \bar N.
\]
Thus, for $t < 1$, there exists $u_t \in C^\infty(\bar N)$ such that $u_t \leq \bar u$ in $\bar N$, $u_t = \varphi$ on $\partial N$ and 
\[
f(\lambda(A_{g_{u_t}})) \geq t\psi(\cdot,u_t) \text{ in } \bar N.
\]

As mentioned above, $\|\ln u_t\|_{C^1(\bar N)}$ is uniformly bounded as $t \rightarrow 1$. Furthermore, as $u_t \leq \bar u$, known interior first derivative estimates \cite{Chen05, GW03-IMRN}, \cite[Theorem 1.10]{Li09-CPAM}, \cite{Wang06} and interior second derivative estimates in \cite{GW03-IMRN}, \cite[Theorem 1.20]{LiLi03} give
\[
\|\ln u_t\|_{C^l(K)} \leq C_{K,l},
\]
for every compact subset $K$ of $N$ and every $l \geq 1$, where $C_{K,l}$ is some constant independent of $t$. Consequently, along a sequence $t_j \rightarrow t$, $\{u_{t_j}\}$ converges in $C^\infty_{loc}(N)$ to some solution $u \in C^\infty(N) \cap C^{0,1}(\bar N)$ of \eqref{Eq:NondegEq}-\eqref{Eq:NdBC}. The proof is complete.
\end{proof}


\subsection{Construction of super-solutions}\label{Sec:SupSolCon}

The following gives a super-solution for Green's functions with a single pole of unit strength. It is clear that a similar construction can be done for any finite number of poles.

\begin{proposition}\label{Prop:SuperSol}
Let $(M,g)$ be an $n$-dimensional smooth compact Riemannian manifold with $n \geq 3$. Let $\Gamma$ satisfy \eqref{G1}, \eqref{G2} and that $\mu_\Gamma^+ > 1$. Assume that $\lambda(A_g) \in \Gamma$ in $M$. Then, for every $p \in M$, there exists a function $\bar u_p \in C^\infty(M \setminus \{p\})$ such that 
\begin{align}
&\lambda\big(A_{g_{\bar u_p}}\big) \in \Gamma \text{ and } \bar u_p > 0 \text{ in } M \setminus \{p\},\label{Eq:PSSP1}\\
&\lim_{x \rightarrow p} d_g(x,p)^{n-2} \bar u_p(x) = 1.\label{Eq:PSSP2}
\end{align}
Furthermore, for every $\mu \in (1,\mu_\Gamma^+] \cap (1,3)$ and $\delta \in (\mu,3)$ and for every sufficiently small $r_1 > 0$, one can arrange, for some $a > 0$, that 
\begin{align}
\bar u_p(x) &= (d_g(x,p)^{-\mu + 1} + a - d_g(x,p)^{-\mu + \delta})^{\frac{n-2}{\mu - 1}}
 \text{ for } 0 < d_g(x,p) < r_1/2,
 \label{Eq:PSSP3}\\
\bar u_p(x) & = 1 \text{ for } d_g(x,p) > r_1.
\label{Eq:PSSP4}
\end{align}
\end{proposition}

\begin{proof} Fix $p \in M$. Let $r(x) := d_g(x,p)$. Fix some $\mu \in (1, \mu_\Gamma^+]$ and $\delta \in (\mu,3)$. Consider, for $a > 1$, the functions
\[
v = v_a = (r^{-\mu + 1} + a - r^{-\mu + \delta})^{\frac{n-2}{\mu - 1}}.
\]
We will show that there exists some $r_2 \in (0,1)$ and $a_0 > 1$ such that
\begin{equation}
\lambda_g\Big(A_{g_{v_a}} - A_g\Big) \in \Gamma \text{ in } \{0 < r < r_2\} \text{ for all } a > a_0,
	\label{Eq:ConBB}
\end{equation}
where $\lambda_g$ signifies that the eigenvalues are computed with respect to $g$. 

We adapt the proof of \cite[Lemma 3.5]{LiNgPoorMan}; the main difference is to allow the possibility that $\mu = \mu_\Gamma^+$. In the sequel, $C$ denotes some positive constant which will always be independent of $a$. Observe that, in local normal coordinates $x^1 = x_1, \ldots, x^n = x_n$ at $p$, the $(0,2)$-Schouten tensor of the metric $v^{\frac{4}{n-2}}g$ satisfies
\[
 A_{g_v} = \chi_1\,g - \chi_2 \frac{x}{r} \otimes \frac{x}{r} +  A_g  + \textrm{err}_1 + \textrm{err}_2,
\]
where $x \otimes x = x_i\,x_j\,dx^i\,dx^j$, 
\begin{align*}
\chi_1 
	&= -\frac{2}{n-2} \frac{v'}{rv} - \frac{2}{(n-2)^2} \frac{|v'|^2}{v^2} \\
	&= \frac{2((\mu-1) + (-\mu + \delta)r^{2})((\mu-1)a - (\delta - 1)\,r^{-\mu + \delta})}{(\mu-1)^2r^{3-\mu}(1 + a\,r^{\mu - 1} - r^{\delta - 1})^2},\\
\chi_2
	&=  \frac{2}{n-2} \frac{1}{v}(v'' - \frac{v'}{r}) - \frac{2n}{(n-2)^2} \frac{|v'|^2}{v^2} \\
	&=  (\mu+1)\chi_1 - \frac{2(\delta - 1)(\delta-\mu)}{(\mu-1)r^{3-\delta}(1 + a\,r^{\mu - 1} - r^{\delta - 1})}
\end{align*}
and
\begin{align*}
|\textrm{err}_1| 
	& \leq C\,r^2\,|\chi_1|,\\
|\textrm{err}_2| 
	&\leq C(r\,v^{-1}|v'| + r^2\,v^{-2}\,|v'|^2) \leq \frac{C}{1 + a r^{\mu - 1}}.
\end{align*}

As $1 < \mu < \delta$, we can assume that $a_0$ is sufficiently large and $r_2$ is sufficiently small such that
\[
\chi_1 \geq \frac{a}{Cr^{3-\mu}(1 + a\,r^{\mu - 1})^2} > 0.
\]
It is important to note that, as $g_{ij}x^i = x^j$, $\lambda_g(\chi_1\,g - \chi_2 \frac{x}{r} \otimes \frac{x}{r}) = (\chi_1 - \chi_2, \chi_1, \ldots \chi_1)$
and so, as $(\chi_1 - \chi_2) > -\mu \chi_1 \geq -\mu_\Gamma^+ \chi_1$ (where we have used $\delta > \mu$), we have
\[
\lambda_g(\chi_1\,g - \chi_2 \frac{x}{r} \otimes \frac{x}{r}) \in \Gamma.
\]
We would like to turn the above relation into a more quantitative form so that it can be used to control the error term.

We have
\begin{equation}
|\textrm{err}_2| \leq Cr^{3-\mu}(a^{-1} + r^{\mu -1})\chi_1.
	\label{Eq:SupSErr}
\end{equation}
For sufficiently large $a_0$ and sufficiently small $r_0$, the right hand side of \eqref{Eq:SupSErr} is smaller than $\chi_1$. Thus, as $\lambda_g(A_g) \in \Gamma$ in $M$ and $M$ is compact, there exists $\nu_0 > 0$ such that
\[
\lambda_g(A_g + \textrm{err}_1 + \textrm{err}_2) \in \Gamma \text{ wherever } \chi_1 < \nu_0.
\]
Thus, by Lemma \ref{Lem:CoConeMatrix}, we have
\[
\lambda_g\Big(A_{g_{v_a}} - A_g\Big) \in \Gamma \text{ in } \{x : 0 < r(x) < r_2, \chi_1(x) < \nu_0\}.
\]

We assume henceforth that $\chi_1 \geq \nu_0$. We have, as $\mu < \delta$,
\[
(1 - \frac{\chi_2}{\chi_1}) + \mu \geq \frac{1}{C}\,r^{\delta - \mu}\,(a^{-1} + r^{\mu - 1}),
\]
which implies in view of the definition of $\mu_\Gamma^+$ and the fact that $\mu \leq \mu_\Gamma^+$ that
\[
\textrm{dist}\Big((1 - \frac{\chi_2}{\chi_1}, 1, \ldots, 1), \RR^n \setminus \Gamma\Big) \geq C\,\min(r^{\delta - \mu}\,(a^{-1} + r^{\mu - 1}),1),
\]
and, as $\chi_1 \geq \nu_0$,
\begin{equation}
\textrm{dist}\Big((\chi_1 - \chi_2, \chi_1, \ldots, \chi_1), \RR^n \setminus \Gamma\Big) \geq C\,\min(r^{\delta - \mu}\,(a^{-1} + r^{\mu - 1}),1)\chi_1.
	\label{Eq:SupSGain}
\end{equation}

The eigenvalues $\tilde\lambda = (\tilde \lambda_1, \ldots, \tilde\lambda_n)$ of $A_{g_v} - A_g$ with respect to the metric $g$ satisfy (cf. \cite[Lemma A.1]{LiNgPoorMan})
\[
|\tilde\lambda_1 - (\chi_1 - \chi_2)| + \sum_{i=2}^n |\tilde\lambda_i - \chi_1|
	\leq C(|\textrm{err}_1| + |\textrm{err}_2|) \leq C[r^{3-\mu}(a^{-1} + r^{\mu -1}) + r^2]\chi_1.
\]
Hence, in view of \eqref{Eq:SupSGain}, we deduce that there is some $r_2 > 0$ and $a_0 > 1$ such that $\tilde \lambda \in \Gamma$ in $\{0 < r < r_2\}$ for $a > a_0$. As $\lambda(A_g) \in \Gamma$, the assertion \eqref{Eq:ConBB} is then readily seen from Lemma \ref{Lem:CoConeMatrix}.

We now turn to the construction of $\bar u_p$. Fix some $\xi \in (0,\frac{1}{10})$. In what follows, the constants $C$ will be also independent of $\xi$. We assume also that $r_1 \in (0,r_2)$ is sufficiently small so that

\begin{equation}
a < v^{\frac{\mu-1}{n-2}} < (1 + \xi)a \text{ and } \frac{|v'| + |v''|}{v} \leq \xi \text{ in } \{r_1/2 < r < r_1\} \text{ for all } a > \frac{C}{\xi}r_1^{-\mu-1}.
	\label{Eq:Conxi}
\end{equation}

Fix some $\varphi \in C_c^\infty(\{r < \frac{4}{5}r_1\})$ such that $\varphi \equiv 1$ in $\{r < \frac{3}{5} r_1\}$ and define
\[
\bar u_p = v_a\,\varphi + a(1- \varphi).
\]
To conclude the proof it suffices to check that, for some sufficiently large $a > 1$,
\begin{equation}
\lambda\big(A_{g_{\bar u_p}}\big) \in \Gamma \text{ in } \{r_1/2 < r < r_1\}.
	\label{Eq:ConAA}
\end{equation}

Using \eqref{Eq:Conxi}, we compute in $\{r_1/2 < r < r_1\}$,
\begin{align*}
\Big|\bar u_p^{-1}\nabla_g \bar u_p - v^{-1}\nabla_g v\Big|_g
	&= \Big|\frac{a(\varphi - 1)}{v \bar u_p} \nabla_g v  + \frac{v - a}{\bar u_p} \nabla_g \varphi\Big|_g\\
	&\leq C\xi ,\\
\Big|\bar u_p^{-1}\nabla_g^2 \bar u_p - v_\vareps^{-1}\nabla_g^2 v_\vareps\Big|_g
	&= \Big|\frac{a(\varphi - 1)}{v\,\bar u_p} \nabla_g^2 v  + \frac{v - a}{\bar u_p} \nabla_g^2 \varphi
	+  \frac{1}{\bar u_p}(d v \otimes d \varphi + d \varphi \otimes d v) \Big|_g\\
	&\leq C\xi.
\end{align*}
Thus, we can write
\[
A_{g_{\bar u}} 
	= (A_{g_v} - A_g) + (A_g + O(\xi)).
\]
We now choose $\xi$ sufficiently small such that $\lambda(A_g + O(\xi)) \in \Gamma$ in $\{r < r_2\}$ and then fix some $a > \max(a_0,\frac{C}{\xi}\,r_1^{-\mu - 1})$ (recall \eqref{Eq:Conxi}). The above computation is then valid, yielding \eqref{Eq:ConAA} as desired.
\end{proof}

\subsection{Existence}\label{ssec:Existence}

Fix $p \in M$ and let $r(x) = d_g(x,p)$. Let $\bar G_p$ be the unique smooth solution of
\begin{align*}
&-\Delta_g \bar G_p + \frac{n-2}{4(n-1)}R_g\,\bar G_p = 0 \text{ and } \bar G_p > 0 \text{ in } M \setminus \{p\},
	\\
&\lim_{x \rightarrow p} r(x)^{n-2} \bar G_p(x) = 1.
\end{align*}
It is well known that such $\bar G_p$ exists and furthermore (cf. \cite{Lee-Parker}),
\begin{equation}
\bar G_p = r^{2-n}(1 + O(r)) \text{ as } r \rightarrow 0.
	\label{Eq:GpAs}
\end{equation}
It should be clear that $\bar G_p = \frac{1}{(n-2) |\SSphere^{n-1}|}G_p$, where $G_p$ is the Green's function for the conformal Laplacian with pole at $p$.

If $\Gamma = \Gamma_1$, we are done. Assume from now on that $\Gamma \neq \Gamma_1$.

Let $f$ be as in Proposition \ref{prop:fConstr}.

Let $\bar u_p$ be as in Proposition \ref{Prop:SuperSol} for some $\mu \in (1,2)$. By \eqref{Eq:PSSP3} and \eqref{Eq:GpAs}, there exists some $r_0 > 0$ such that $\bar u_p > \bar G_p$ in $\{ 0 < r < r_0\}$. On the other hand, by \eqref{Eq:PSSP1}, 
\[
-\Delta \bar u_p + \frac{n-2}{4(n-1)}R_g\,\bar u_p \geq 0 \text{ in } M \setminus \{p\}.
\]
Hence, by the maximum principle,
\begin{equation}
\bar u_p \geq \bar G_p \text{ in } M \setminus \{p\}.
	\label{Eq:up>Gp}
\end{equation}

For small $\delta > 0$, let
\[
c_\delta = \min_{M \setminus B_\delta(p)} f\big(\lambda(A_{g_{\bar u}})\big) > 0.
\]
(Here we have used the smoothness of $\bar u_p$ to establish the positivity of $c_\delta$.) By Theorem \ref{Thm:GuanGen}, for every $c \in (0,c_\delta)$, there exists a function $u_{\delta,c} \in C^{0,1}(M \setminus B_\delta(p)) \cap C^\infty(M \setminus \bar B_\delta(p))$ satisfying
\begin{align}
&f\big(\lambda(A_{g_{u_{\delta,c}}})\big) = c \text{ and } u_{\delta,c} > 0 \text{ in } M \setminus B_\delta(p),\label{Eq:AppXeq}\\
&u_{\delta,c} = \bar u_p \text{ on } \partial B_\delta(p),\label{Eq:AppXBC}
\end{align}
Furthermore, $\{\ln u_{\delta,c}\}_{c \in (0,c_\delta)}$ is uniformly bounded in $C^{0,1}(M \setminus B_\delta(p))$ and $C^\infty_{loc}(M \setminus \bar B_\delta(p))$. It follows that, along a sequence $c_j \rightarrow 0$, $\{u_{\delta,c_j}\}$ converges in $C^2_{loc}(M \setminus \bar B_\delta(p))$ to some smooth functions $u_\delta \in C^{0,1}(M \setminus B_\delta(p)) \cap C^\infty(M \setminus \bar B_\delta(p))$ 
satisfying
\begin{align}
&\lambda\big(A_{g_{u_{\delta}}}\big) \in \partial\Gamma \text{ and } u_{\delta} > 0 \text{ in } M \setminus B_\delta(p),\label{Eq:AppXeqDeg}\\
&u_{\delta} = \bar u_p \text{ on } \partial B_\delta(p),\label{Eq:AppXBCDeg}
\end{align}

Using \eqref{Eq:PSSP1}, \eqref{Eq:AppXeqDeg} and the maximum principle, we see that 
\[
u_\delta \leq \bar u_p \text{ in } M \setminus B_\delta(p).
\]
Hence, for each compact subset $K$ of $M \setminus \{p\}$, there exist constants $C_K > 0 $ and $0 < \bar c_{\delta, K} < c_\delta$ such that 
\[
u_{\delta, c} \leq C_K \text{ provided } c < \bar c_{\delta, K}.
\]
By known first and second derivative estimates, for every compact subset $K'$ of $\mathring{K}$, there holds
\[
\|u_{\delta, c}\|_{C^2(K')} \leq C_{K,K'} \text{ for all } c < \bar c_{\delta, K},
\]
where $C_{K,K'}$ is independent of $\delta$. Sending $c$ to zero, we obtain that
\[
\|u_{\delta}\|_{C^2(K')} \leq C_{K,K'}.
\]
In other words, the family $\{u_\delta\}$ is bounded in $C^2_{loc}(M \setminus \{p\})$. Hence, there is some $\delta_j \rightarrow 0$ such that $\{u_{\delta_j}\}$ converges in $C^{1,\alpha}_{loc}(M \setminus \{p\})$ for any $\alpha \in (0,1)$ to some $u \in C^{1,1}_{loc}(M \setminus \{p\})$.

As $u_\delta \leq \bar u_p$, we have $u \leq \bar u_p$ in $M \setminus \{p\}$. On the other hand, by \eqref{Eq:AppXeq},
\[
-\Delta_g u_{\delta,c} + \frac{n-2}{4(n-1)}R_g\,u_{\delta,c} \geq 0 \text{ in } M \setminus B_\delta(p).
\]
In view of \eqref{Eq:up>Gp}, \eqref{Eq:AppXBC} and the maximum principle, we thus have $u_{\delta, c} \geq \bar G_p$ in $M \setminus B_\delta(p)$. It follows that $u \geq \bar G_p$ in $M \setminus \{p\}$. On one hand, this implies \eqref{Eq:DegAs}. On the other hand, this implies $u > 0$ in $M \setminus \{p\}$, and so by \eqref{Eq:AppXeqDeg} and the convergence of $\{u_{\delta_j}\}$ to $u$, we obtain \eqref{Eq:DegEq1}. We have thus proved the existence of a solution to \eqref{Eq:DegEq1}-\eqref{Eq:DegAs}.

\begin{remark}\label{Rem:6.1}
By construction, we have $\bar G_p \leq u \leq \bar u_p$. Hence,  for any $\mu \in (1, \mu_\Gamma^+] \cap (1,2)$, one has
\[
0 \leq \liminf_{x \rightarrow p} \frac{r(x)^{n-2}\,u(x) - 1}{r(x)^{\mu - 1}} \leq \limsup_{x \rightarrow p} \frac{r(x)^{n-2}\,u(x) - 1}{r(x)^{\mu - 1}} < \infty.
\]
(If there are multiple poles with multiple strengths, we have
\[
0 \leq \liminf_{x \rightarrow p_i} \frac{d_g(x,p_i)^{n-2}\,u(x) - c_i}{d_g(x,p_i)^{\mu - 1}} \leq \limsup_{x \rightarrow p_i} \frac{d_g(x,p_i)^{n-2}\,u(x) - c_i}{d_g(x,p_i)^{\mu - 1}} < \infty.)
\]
When $\mu = \mu_\Gamma^+ < 2$, this is in a sense sharp. See \cite[Theorem 1.2]{LiNgBocher}, where it is shown that if $\lambda(A_{U^{\frac{4}{n-2}}g_{Euc}}) \in \partial\Gamma$ on a punctured ball of the flat space $(\RR^n,g_{Euc})$ and if $\mu_\Gamma^+ > 1$ and $(1, 0, \ldots, 0) \in \partial\Gamma$, then $U$ can be expressed in the form
\[
U(x) = \Big(c\,|x|^{-\mu_\Gamma^{+} + 1} + \mathring{w}(x)\Big)^{\frac{n-2}{\mu_\Gamma^+ - 1}}
\]
for some non-negative bounded function $\mathring{w}$ which is either positive or identically zero.
\end{remark}

\subsection{Uniqueness}\label{ssec:Uniqueness}

In this subsection, we prove that \eqref{Eq:DegEq1}-\eqref{Eq:DegAs}  has a unique continuous viscosity solution. Let $u \in C^{1,1}_{loc}(M \setminus \{p\})$ be the solution to \eqref{Eq:DegEq1}-\eqref{Eq:DegAs} which was constructed in subsection \ref{ssec:Existence}. Assume that $v \in C^{0}_{loc}(M \setminus \{p\})$ is also a solution to \eqref{Eq:DegEq1}-\eqref{Eq:DegAs}.

\medskip
\noindent{\it Step 1.} We show that $v \leq u$. To this end, we show that $\theta v \leq u$ for all $\theta \in (0,1)$.

By construction, there exist sequences $r_j \rightarrow 0$, $\vareps_j \rightarrow 0$ and $\{u_j\} \subset C^\infty(M \setminus B_{r_j}(p))$ such that $\{u_j\}$ converges to $u$ in $C^{1,\alpha}_{loc}(M \setminus S)$ and
\begin{align}
&f\big(\lambda(A_{g_{u_j}})\big) = \vareps_j \text{ and }  u_j > 0 \text{ in } M \setminus B_{r_j}(p),
	\label{Eq:fuj=ej}\\
&\lim_{j \rightarrow \infty} r_j^{n-2} \sup_{\partial B_{r_j}(p)} u_j = \lim_{j \rightarrow \infty} r_j^{n-2} \inf_{\partial B_{r_j}(p)} u_j = 1.
	\label{Eq:fujBC}
\end{align}

Clearly, by \eqref{Eq:DegAs} and \eqref{Eq:fujBC}, for sufficiently large $j$, $\theta v < u_j$ on $\partial B_{r_j}(p)$. We claim that $\theta v \leq u_j$ in $M \setminus \bar B_{r_j}(p)$. Indeed, if this is not true, there is some $\alpha \in (0,1)$ and $q \in M \setminus \bar B_{r_j}(p)$ such that $\alpha\theta v \leq u_j$ in $M \setminus \bar B_{r_j}(p)$ and $\alpha\theta v(q) \leq u_j(q)$. As $\lambda\big(A_{g_{\alpha\theta v}}\big) \in \partial \Gamma$ and $u_j$ is smooth, it follows that
\[
\lambda(A_{g_{u_j}}) \in \RR^n \setminus \Gamma,
\]
which contradicts \eqref{Eq:fuj=ej}. We have thus shown that $\theta v \leq u_j$ in $M \setminus \bar B_{r_j}(p)$. Sending $j \rightarrow \infty$ and then $\theta \rightarrow 1$, we arrive at $v \leq u$ in $M \setminus \{p\}$.

\medskip
\noindent{\it Step 2.} We show that $v \geq u$. Similar to the previous step, we in fact show that $v \geq \theta u$ for all $\theta \in (0,1)$.

Clearly, there exists some $r_0 > 0$ such that 
\[
v > \theta u \text{ in } \bar B_{r_0}(p) \setminus \{p\}.
\]

Let $u_j$ be as before in Step 1. It is more convenient to work with $w = u^{-\frac{2}{n-2}}$, $w_j = u_j^{-\frac{2}{n-2}}$, and $\xi = v^{-\frac{2}{n-2}}$. We then have
\begin{align}
&\lambda_g(A_w), \lambda_g(A_\xi) \in \partial\Gamma \text{ in } M \setminus \{p\},\label{Eq:XYZ1}\\
&f(\lambda_g(A_{w_j})) = \vareps_j w_j^{-1}\text{ in } M \setminus B_{r_j}(p),\label{Eq:XYZ2}
\end{align}
where here and below $\lambda_g$ signifies that the eigenvalues are computed with respect to the metric $g$ and
\[
A_\psi = \nabla^2 \psi - \frac{1}{2\psi}|d \psi|^2_g\,g + \psi\,A_g.
\]

As $\{w_j\}$ converges in $C^0_{loc}(M \setminus \{p\})$ to $w$, which is positive on $M \setminus \{p\}$, there exists some $\bar \mu > 0$ such that, for all sufficiently large $j$,
\[
w_j > 2\bar\mu \text{ in } M \setminus B_{r_0}(p).
\]
Fix some $\mu \in (0,\bar \mu)$ for the moment. We have
\[
A_{w_j} 
	= A_{w_j - \mu} +  \frac{\mu}{2w_j(w_j - \mu)} |dw_j|_g^2\,g + \mu\,A_g.
\]
As $M$ is compact and $\lambda_g(A_g) \in \Gamma$, there is some $\delta > 0$ such that 
\[
\lambda_g(A_g - 2\delta g) \in \Gamma \text{ in } M.
\]
We now write, 
\[
A_{w_j} - \delta g = A_{w_j - \mu} + (A_g - 2\delta g) +  \Big(\delta + \frac{\mu}{2w_j(w_j - \mu)} |dw_j|_g^2\Big)g.
\]
On the other hand, by \eqref{Eq:XYZ2} and the fact that $\{w_j\}$ is uniformly bounded in $C^2(M \setminus B_{r_0}(p))$, $\lambda_g(A_{w_j} - \delta g) \in \RR^n \setminus \bar\Gamma$ in $M\setminus B_{r_0}(p)$ for all sufficiently large $j$. Invoking Lemma \ref{Lem:CoConeMatrix} again, we thus have
\begin{equation}
\lambda(A_{w_j - \mu}) \in \RR^n \setminus \bar \Gamma \text{ in } M \setminus B_{r_0}(p)
	\text{ for all sufficiently large $j$}.
		\label{Eq:wjPushout}
\end{equation}

Using \eqref{Eq:wjPushout}, we can argue as in Step 1 to show that, for all sufficiently large $j$,
\[
\xi \leq \theta^{-\frac{n-2}{2}}(w_j - \mu) \text{ in } M \setminus B_{r_0}(p).
\]
Sending $j \rightarrow \infty$ and then $\mu \rightarrow 0$, we obtain that
\[
v \geq \theta u \text{ in } M \setminus B_{r_0}(p).
\]
Recalling the definition of $r_0$, we conclude that $v \geq \theta u$ in $M \setminus \{p\}$, which upon letting $\theta \rightarrow 1$ yields $v \geq u$ in $M \setminus \{p\}$.

Combining Step 1 and Step 2, we conclude that $v \equiv u$, i.e. the solution to  \eqref{Eq:DegEq1}-\eqref{Eq:DegAs} is unique. This completes the proof of Theorem \ref{thm:MainThm}.\hfill$\Box$



\section{Green's functions as solutions to nonlinear equations with Dirac delta measures on the right hand side}
\label{Sec:delta}

In this section, we illustrate that Green's functions may show up as suitable rescaled limits for certain blow-up solutions to the nonlinear Yamabe problem
\begin{equation}
f\big(\lambda(A_{g_u})\big) = 1 \text{ and } u > 0.
	\label{Eq:04VI19-NE}
\end{equation}
More general scenarios of blow-up will be considered elsewhere.

Let $(M,g)$ be a compact Riemannian manifold and $i(M,g)$ its injectivity radius. Suppose for some $0 < r_0 < i(M,g)$ that $\{u_i\}$ is a sequence of smooth solutions to \eqref{Eq:04VI19-NE} on some balls $B_g(p_i,r_0)$ in $M$ such that $u_i(p_i) = \max_{B_g(p_i,r_0)} u_i \rightarrow \infty$ and $p_i \rightarrow p_\infty$ as $i \rightarrow \infty$. We say that $\{u_i\}$ has an isolated blow-up point if 
\begin{enumerate}[({H}1)]
\item there exists $C > 0$ independent of $i$ such that $ d_g(\cdot, p_i)^{\frac{n-2}{2}} u_i \leq C$ in $B_{g}(p_i,r_0)$.
\end{enumerate}
We say that $\{u_i\}$ has tame geometry in $B_g(p_i,r_0)$ if 
\begin{enumerate}[({H}1)]
\setcounter{enumi}{1}
\item there exist $C > 0$ and $\theta \in [0,1)$ independent of $i$ such that 
\begin{equation}
|\Ric_{g_{u_i}}|_{g_{u_i}} \leq C \max(1, u_i(p_i)^{\frac{4\theta}{n-2}} d_g(\cdot, p_i)^{2\theta}) \text{ in } B_{g}(p_i,r_0).
	\label{Eq:13VIII20-C1}
\end{equation}
\end{enumerate}
Note that, by \cite[Proposition B.1]{LiNgBocher}, for $\Gamma = \Gamma_k$ with $2 \leq k \leq n$, \eqref{Eq:13VIII20-C1} can be replaced by
\[
R_{g_{u_i}} \leq C \max(1, u_i(p_i)^{\frac{4\theta}{n-2}} d_g(\cdot, p_i)^{2\theta}) \text{ in } B_{g}(p_i,r_0).
	\tag{\ref{Eq:13VIII20-C1}'}
\]

When $\theta = 0$ in $(H2)$, we say that $\{u_i\}$ has bounded geometry. It should be noted that, by \cite{LiNgPoorMan}, when $(M,g)$ is not conformal to the standard sphere, equation \eqref{Eq:04VI19-NE} on $M$ has no blow-up sequence of solutions with bounded geometry on the whole of $M$.

It should also be noted that, under $(H1)$, it is easy to show (in view of estimate \eqref{Eq:08VII19-A1} and Lemma \ref{Lem:08VI19-LowerSharpDecay}) that estimate \eqref{Eq:13VIII20-C1} holds with $\theta = 1$, i.e.
\[
|\Ric_{g_{u_i}}|_{g_{u_i}} \leq C \max(1, u_i(p_i)^{\frac{4}{n-2}} d_g(\cdot, p_i)^{2}) \text{ in } B_{g}(p_i,r_0/2).
\]
It is clear from the above that, under $(H1)$, if $(H2)$ holds for some $\theta = \theta_0$, then it holds for all $\theta \in (\theta_0,1)$, after a shrinking $r_0$. We do not know yet whether $(H1)$ implies $(H2)$ in general.

\begin{theorem}
\label{Thm:deltaRHS-IsoTGeoBlUp}
Let $(M,g)$ be an $n$-dimensional smooth compact Riemannian manifold with $n \geq 3$. Suppose that $(f,\Gamma)$ satisfies \eqref{G1}-\eqref{f4} and that $\muGp > 1$. Suppose that $\{u_i\}$ is a sequence of solutions to \eqref{Eq:04VI19-NE} on some balls $B_g(p_i,r_0) \subset M$ with $0 < r_0 < i(M,g)$ independent of $i$ which has an isolated blow-up point and has tame geometry (i.e. $(H1)$ and $(H2)$ hold). Then, upon extracting a subsequence, $\tilde u_i := u_i(p_i) u_i$ converges in $C^{1,\alpha}_{\rm loc}(B_{g}(p_\infty,r_0/2)\setminus \{p_\infty\})$ for every $0 < \alpha < 1$ to a positive function $\tilde u_\infty \in C^{\infty}(B_{g}(p_\infty,r_0/2)\setminus \{p_\infty\})$ satisfying
\begin{align*}
&A_{g_{\tilde u_\infty}} \equiv 0 \text{ in } B_{g}(p_\infty,r_0/2) \setminus \{p_\infty\},\\
&\lim_{x \rightarrow p_\infty} d(x, p_\infty)^{n-2}\tilde u_\infty(x) \in (0,\infty).
\end{align*}
Furthermore, if $(f,\Gamma) = (\sigma_k^{1/k}, \Gamma_k)$ for some $1 \leq k < \frac{n}{2}$, then $\tilde u_i^{\frac{n+2k}{n-2}} \sigma_k\big(\lambda(A_{g_{\tilde u_i}})\big)$ weakly* converges in $B(p_\infty,r_0/2)$ in the space of measures to a Dirac measure $m_{n,k} \delta_{p_\infty}$ with an explicit $m_{n,k} > 0$.
\end{theorem}

The conclusion of the above theorem holds if we replace the right hand side of \eqref{Eq:04VI19-NE} by a smooth positive function $\eta(x)$, in which case the limit measure changes to $\eta(p_\infty)^{\frac{n-2}{2}}m_{n,k} \delta_{p_\infty}$.

The rest of the section contains 2 subsections. The proof of Theorem \ref{Thm:deltaRHS-IsoTGeoBlUp} is given in Subsection \ref{SSec:IsoTBlUpSeq}. We first show that $(H1)$ and $(H2)$ rule out a phenomenon usually known as bubbles on top of bubbles; see Lemma \ref{Lem:07VIII19-Iso+TGeo->Sim}. Using a suitable barrier construction, we then show a sub-optimal upper bound for $u_i$ (see \eqref{Eq:07VIII19-E1}) which is sufficient to establish the weak* convergence of $\tilde u_i^{\frac{n+4}{n-2}} \sigma_2(\lambda(A_{g_{\tilde u_i}}))$ and to identify its limit; see Lemma \ref{Lem:08VI19-AlmostSharpDecay} and Corollary \ref{Cor:sigmakL1Bnd}. Exploiting further condition $(H2)$, we then derive a sharper upper bound for $u_i$ in Lemma \ref{Lem:08VI19-UpperSharpDecay}, and deduce the convergence of $\tilde u_i$, which concludes the proof. In Subsection \ref{SSec:DivId}, we use the divergence structure associated with the $\sigma_k$ operator to prove a compensated compactness type result for the $\sigma_k$ equation (see Proposition \ref{Prop:CompCptness1}). This is not directly related to the proof of Theorem \ref{Thm:deltaRHS-IsoTGeoBlUp} but may be relevant in the study of Green's functions.

\subsection{Isolated blow-up sequences with tame geometry}\label{SSec:IsoTBlUpSeq}

Let $\{u_i\}$ be a sequence of smooth local solutions to the nonlinear Yamabe equation \eqref{Eq:04VI19-NE} 
\[
f\big(\lambda(A_{g_{u_i}})\big) = 1 \text{ and } u_i > 0 \text{ on some ball } B_g(p_i,r_0)
\]
with $0 < r_0 < i(M,g)$. We suppose that $\{u_i\}$ has an isolated blow-up point and has tame geometry, i.e. we have that $u_i(p_i) = \max_{B_g(p_i,r_0)} u_i \rightarrow \infty$, $p_i \rightarrow p_\infty$, and that conditions $(H1)$ and $(H2)$ hold.

We aim to show that $\tilde u_i = u_i(p_i) u_i$ converges to a solution $\tilde u_\infty$ of $\lambda(A_{g_{\tilde u_\infty}}) \in \partial \Gamma$ (in fact $A_{g_{\tilde u_\infty}} \equiv 0$) with $u_\infty(x) = c(1 + o(1)) d_g(x,p_\infty)^{-(n-2)}$ near $p_\infty$ for some constant $c \in (0,\infty)$ and, when $(f,\Gamma) = (\sigma_k^{1/k},\Gamma_k)$, to identify the weak* limit of the sequence $\tilde u_i^{\frac{n+2k}{n-2}} \sigma_k(\lambda(A_{g_{\tilde u_i}}))$.

\subsubsection{Preliminary analysis}

We start with some well-known facts. By local gradient and second derivative estimates (\cite{Chen05, GW03-IMRN,Li09-CPAM,Wang06,LiLi03}), we have
\begin{equation}
|\nabla^\ell \ln u_i(x)| \leq Cd_g(x, p_i)^{-\ell} \text{ in } B_g(p_\infty,3r_0/4) \setminus \{p_i\} \text{ for } \ell = 1,2.
	\label{Eq:08VII19-A1}
\end{equation}

For $x \in \RR^n$ and $\lambda > 0$, let
\[
U_\lambda(x) = \varkappa \Big(\frac{\lambda}{1 + \lambda^2 |x|^2}\Big)^{\frac{n-2}{2}},
\]
where $|\cdot|$ denotes the Euclidean norm and $\varkappa = \varkappa(f,\Gamma)$ is a (normalizing) positive constant so that
\[
f\big(\lambda(A_{\mathring{g}_{U_\lambda}})\big) = 1 \text{ on } \RR^n \text{ for all } \lambda > 0,
\]
where $\mathring{g}$ denotes the Euclidean metric on $\RR^n$.
 
Define a map $\Phi_i: \RR^n \approx T_{p_i}(M,g) \rightarrow M$ by 
\[ 
\Phi_i(x) = \exp_{p_i} \frac{\varkappa^{\frac{2}{n-2}}\,x}{u_i(p_i)^{\frac{2}{n-2}}},
\] 
and let  
\[ 
\hat u_i(x) = \varkappa u_i(p_i)^{-1}\,u_i \circ \Phi_i(x), \qquad x  \in \RR^n. 
\] 
Then $\hat u_i$ satisfies 
\begin{equation}
\sigma_k(\lambda(A_{(\hat g_i)_{\hat u_i}})) = 1  \text{ in } \{|x| < \delta_0\,\varkappa^{-\frac{2}{n-2}}\,u_i(p_i)^{\frac{2}{n-2}}\}, 
	\label{Eq:05VI19-A2}
\end{equation} 
where $\hat g_i := \varkappa^{-\frac{4}{n-2}} u_i(p_i)^{-\frac{4}{n-2}} \Phi_i^* g$ and $\delta_0$ is the injectivity radius of $(M,g)$. It is clear that $\hat g_i \rightarrow \mathring{g}$ in $C^3_{\rm loc}(\RR^n)$. Furthermore $\hat u_i(0) = \varkappa$ and $\hat u_i \leq \varkappa$ in $\{|x| < \delta_0\,\varkappa^{-\frac{2}{n-2}}\,u_i(p_i)^{\frac{2}{n-2}}\}$. By known local first and second derivative estimates, it follows that $\{\ln \hat u_i\}$ is uniformly bounded in $C^2$ on any compact subset of $\RR^n$. By Evans-Krylov's theorem and the Schauder theory, $\{\hat u_i\}$ is uniformly bounded in $C^3$ on any compact subset of $\RR^n$ and converges, along a subsequence, in $C^{2,\alpha}_{\rm loc}(\RR^n)$ to some positive $\hat u_* \in C^2(\RR^n)$ which satisfies $\varkappa = \hat u_*(0) = \max \hat u_*$ and
\[ 
\sigma_k(\lambda(A_{\mathring{g}_{\hat u_*}})) = 1  \text{ on } \RR^n.
\] 
By the Liouville theorem \cite[Theorem 1.3]{LiLi05}, we have $\hat u_* = U_{1}$. In particular, passing to another subsequence if necessary, we have for an arbitrarily fixed $N > n$ that $i \ll u_i(p_i)^{\frac{2}{n-2}}$ and
\begin{equation} 
\lim_{i \rightarrow \infty} i^{N} \|\hat u_i - U_1\|_{C^2(\bar B_i)} = 0.
	\label{Eq:05VI19-A3} 
\end{equation} 

\begin{lemma}\label{Lem:08VI19-LowerSharpDecay}
Under the assumptions of Theorem \ref{Thm:deltaRHS-IsoTGeoBlUp} except for $(H2)$, there exists $C > 1$ (independent of $i$) such that, after passing to a subsequence,
\[
u_i(x) \geq \frac{1}{C} u_i(p_i)^{-1} d_g(x,p_i)^{-(n-2)} \text{ in } \{r_0 \geq d_g(x,p_i) \geq \varkappa^{\frac{2}{n-2}}  u_i(p_i)^{-\frac{2}{n-2}}\}.
\]
\end{lemma}

\begin{proof} In the sequel, $C$ denotes some positive constant which will always be independent of $i$.

Let $L_g = \Delta_g - \frac{n-2}{4(n-1)}R_g$ denote the conformal Laplacian of $g$. We have that $L_g u_i \geq 0$. A calculation shows that there exist large $K$ and small $\delta$ such that, for every $p$ near $p_\infty$, the function $\tilde G_p(x) := d_g(p,x)^{2-n} - K d_g(p,x)^{\frac{5}{2}-n} - (\delta^{2-n} -  K\delta^{\frac{5}{2}-n})$ satisfies (see e.g. \cite[Lemma 3.3]{LiNgPoorMan})
\[
L_g \tilde G_p \geq 0 \text{ in } B(p,\delta) \setminus \{p\}.
\]
Now note that, by \eqref{Eq:05VI19-A3} and with $r_i = \varkappa^{\frac{2}{n-2}}u_i(p_i)^{-\frac{2}{n-2}}$, we have for large $i$ that
\begin{align*}
u_i(x) \geq \frac{1}{C} u_i(p_i)^{-1} \tilde G_{p_i}(x) \text{ on } \partial B(p_i,r_i).
\end{align*}
Clearly $u_i(x) \geq 0 = \frac{1}{C} u_i(p_i)^{-1} \tilde G_{p_i}(x)$ on $\partial B(p_i,\delta)$. An application of the maximum principle then shows that
\[
u_i(x) \geq  \frac{1}{C} u_i(p_i)^{-1} \tilde G_{p_i}(x) \geq \frac{1}{C} u_i(p_i)^{-1} d_g(p_i,x)^{2-n} \text{ in } B(p_i,\delta) \setminus B(p_i,r_i).
\]
The conclusion follows from the gradient estimate \eqref{Eq:08VII19-A1}.
\end{proof}

\subsubsection{Simplicity of blow-up sequences}

In this subsection, we show that if $\{u_i\}$ has only one isolated blow-up point and has tame geometry, then $\{u_i\}$ is simple in the sense that
\begin{enumerate}[({H}1)]
\setcounter{enumi}{2}
\item there exists $r_0' > 0$ independent of $i$ such that the functions
\[
r \mapsto \frac{r^{\frac{n-2}{2}}}{|\partial B_g(p_i, r)|_g}\int_{\partial B_g(p_i, r)} u_i(x)dS_g(x)
\]
are non-increasing in $(2\varkappa^{\frac{2}{n-2}}u_i(p_i)^{-\frac{2}{n-2}},r_0')$.
\end{enumerate}

\begin{lemma}\label{Lem:07VIII19-Iso+TGeo->Sim}
Under the assumptions of Theorem \ref{Thm:deltaRHS-IsoTGeoBlUp}, the sequence $\{u_i\}$ is simple, namely $(H3)$ holds.
\end{lemma}

The proof is by contradiction. We suppose that the sequence $\{u_i\}$ is not simple and rescale it to a situation in which simplicity holds and appeal to the following result.

\begin{lemma}\label{Lem:08VI19-AlmostSharpDecay}
Let $(M,g)$ be an $n$-dimensional smooth compact Riemannian manifold with $n \geq 3$. Suppose that $(f,\Gamma)$ satisfies \eqref{G1}-\eqref{G2}, \eqref{f1}-\eqref{f4} and that $\muGp > 1$. Suppose that $\{u_i\}$ is a sequence of solutions to \eqref{Eq:04VI19-NE} on some balls $B_g(p_i,r_0) \subset M$ with $0 < r_0 < i(M,g)$ independent of $i$ which has an isolated simple blow-up point, i.e. $(H1)$ and $(H3)$ hold. Then, for any $0 < \theta < 1$, there exists $C > 0$ (independent of $i$) such that, after passing to a subsequence,
\begin{equation}
u_i(x) \leq C u_i(p_i)^{-\theta} d_g(x,p_i)^{-\frac{(1 + \theta)(n-2)}{2}} \text{ in } \{ r_0/2 \geq d_g(x,p_i) \geq \varkappa^{\frac{2}{n-2}} i u_i(p_i)^{-\frac{2}{n-2}}\}.
	\label{Eq:07VIII19-E1}
\end{equation}
\end{lemma}

An immediate consequence of the above result for $\theta \in (\frac{n-2k}{n+2k},1)$ is that
\[
u_i(p_i)^{\frac{n-2k}{n-2}} \int_{\{r_0/2 \geq d_g(x,p_i) \geq \varkappa^{\frac{2}{n-2}} i u_i(p_i)^{-\frac{2}{n-2}}\}} 
	u_i(x)^{\frac{n+2k}{n-2}} \,dv_g
	\leq 
	 C i^{\frac{1-\theta}{2}(n+2k)  - 2k} \rightarrow 0.
\]
This together with \eqref{Eq:05VI19-A3} gives:

\begin{corollary}\label{Cor:sigmakL1Bnd}
Under the assumptions of Theorem \ref{Thm:deltaRHS-IsoTGeoBlUp} and with $(f,\Gamma) = (\sigma_k^{1/k}, \Gamma_k)$, we have for any fixed $r < r_0/2$ that
\begin{align*}
\int_{B(p_i,r)} \tilde u_i^{\frac{n+2k}{n-2}}\sigma_k(A_{g_{\tilde u_i}})\,dv_g 
	&= u_i(p_i)^{\frac{n-2k}{n-2}} \int_{B(p_i,r)} u_i^{\frac{n+2k}{n-2}}\,dv_g\\
 	&\rightarrow U_1(0)^{\frac{n-2k}{n-2}} \int_{\RR^n} U_1^{\frac{n+2k}{n-2}}\,dx.
\end{align*}
\end{corollary}

To prove Lemma \ref{Lem:08VI19-AlmostSharpDecay} before that of Lemma \ref{Lem:07VIII19-Iso+TGeo->Sim}, we will need the following lemma.

\begin{lemma}\label{Lem:ConeSuperHar}
Let $(M,g)$ be a smooth compact Riemannian manifold of dimension $n \geq 3$. Assume that $(f,\Gamma)$ satisfies \eqref{G1}-\eqref{G2}, \eqref{f1}-\eqref{f4} and that $\muGp > 1$. For $q \in (0,n-2)$, there exist some $r_1 > 0$ and $C > 1$ such that for every $p \in M$ and $a, b > 0$, the function
\[
\varphi(x) = a r(x)^{-q} + b r(x)^{-(n-2-q)} \text{ in } x \in B_g(p,r_1) \text{ where } r(x) = d_g(x,p)
\]
satisfies
\[
f(\lambda(A_{g_{\varphi}})) \geq \frac{1}{Cr^2} \varphi^{-\frac{4}{n-2}}   \text{ in } B(p,r_1) \setminus \{p\}.
\]
\end{lemma}

\begin{proof} In the sequel, $C$ denotes some positive constant which will always be independent of $a, b$. Observe that, in local normal coordinates $x^1 = x_1, \ldots, x^n = x_n$ at $p$, the Schouten tensor of the metric $g_\varphi$ satisfies
\[
 A_{g_v} = \chi_1\,g - \chi_2 \frac{x}{r} \otimes \frac{x}{r} +  A_g  + \textrm{err}_1 + \textrm{err}_2,
\]
where $x \otimes x = x_i\,x_j\,dx^i\,dx^j$, 
\begin{align*}
\chi_1 
	&= -\frac{2}{n-2} \frac{\varphi'}{r\varphi} - \frac{2}{(n-2)^2} \frac{|\varphi'|^2}{\varphi^2} \\
	&= \frac{2}{(n-2)^2} \frac{(aq r^{-q} + b(n-2-q)r^{-(n-2-q)})(a(n-2-q)r^{-q} + bqr^{-(n-2-q)}))}{r^2\varphi^2} \\
	&\in ( \frac{1}{Cr^2}, \frac{C}{r^2}),\\
\chi_2
	&=  \frac{2}{n-2} \frac{1}{\varphi}(\varphi'' - \frac{\varphi'}{r}) - \frac{2n}{(n-2)^2} \frac{|\varphi'|^2}{\varphi^2} \\
	&=  2\chi_1 - \frac{2ab(n-2-2q)^2}{(n-2)r^n\varphi^2}
\end{align*}
and
\begin{align*}
|\textrm{err}_1| 
	& \leq C,\\
|\textrm{err}_2| 
	&\leq C(r\,v^{-1}|v'| + r^2\,v^{-2}\,|v'|^2) \leq C.
\end{align*}

It follows that the eigenvalues $\lambda = \lambda(A_{g_\varphi})$ (with respect to $g_\varphi$) satisfy
\[
|\lambda_1 - \varphi^{-\frac{4}{n-2}}(\chi_1 - \chi_2)| + \sum_{i = 2}^n |\lambda_i - \varphi^{-\frac{4}{n-2}} \chi_1| 
	\leq C \varphi^{-\frac{4}{n-2}}.
\]
Noting that, as $\muGp > 1$, $(-1, 1, \ldots, 1) \in \Gamma$. It thus follows, for sufficiently small $r_1$, that $\lambda(A_{g_\varphi}) \in \Gamma$ in $\{0 < r < r_1\}$ and
\begin{align*}
f(\lambda(A_{g_\varphi})) 
	&= \varphi^{-\frac{4}{n-2}}\chi_1 f(\chi_1^{-1}(\chi_1  - \chi_2),1, \ldots, 1) + O( \varphi^{-\frac{4}{n-2}})\\
	&= \varphi^{-\frac{4}{n-2}}\chi_1 f(-1,1, \ldots, 1) + O( \varphi^{-\frac{4}{n-2}})\\
	&\geq \frac{1}{Cr^2} \varphi^{-\frac{4}{n-2}} \text{ in } \{0 < r < r_1\},
\end{align*}
which concludes the argument.
\end{proof}

\begin{proof}[Proof of Lemma \ref{Lem:08VI19-AlmostSharpDecay}] In the sequel, $C$ denotes some positive constant which will always be independent of $i$.

Let $r_i = \varkappa^{\frac{2}{n-2}} u_i(p_i)^{-\frac{2}{n-2}}$. By \eqref{Eq:05VI19-A3}, we have
\begin{align*}
u_i(x) \leq C u_i(p_i)^{-1} d_g(x,p_i)^{-(n-2)} \text{ on } \partial B(p_i,i r_i).
\end{align*}
Thus, by isolated simplicity and the gradient estimate \eqref{Eq:08VII19-A1},
\begin{equation}
d_g(x,p_i)^{\frac{n-2}{2}} u_i(x) \leq C u_i(p_i)^{-1} r_i^{-\frac{n-2}{2}} = C i^{-\frac{n-2}{2}}  \text{ in } \{ir_i \leq d_g(x,p_i) \leq r_0\}.
	\label{Eq:08VII19-A2}
\end{equation}
It thus follows, for some constant $C_0 > 0$, that
\begin{equation}
f(\lambda(A_{g_{u_i}})) = 1 \leq \frac{C_0}{i^{2} r^{2}} u_i^{-\frac{4}{n-2}}
\text{ in } \{i r_i \leq d_g(x,p_i) \leq r_0\}.
	\label{Eq:08VII19-B1}
\end{equation}

Let $q =  \frac{(1 - \theta)(n-2)}{2}$. By Lemma \ref{Lem:ConeSuperHar}, for all $a_i,b_i > 0$  the functions
\[
\varphi_i(x) = a_i d_g(x,p_i)^{-q} + b_i d_g(x,p_i)^{-(n-2-q)}
\]
satisfy for some sufficiently small $r_1 \in(0,r_0)$ that
\begin{equation}
 f(\lambda(A_{g_{\varphi_i}})) \geq \frac{1}{C r^{2}} \varphi_i^{-\frac{4}{n-2}} \geq \frac{C_0}{i^{2} r^{2}} \varphi_i^{-\frac{4}{n-2}} \text{ in } \{0 \leq d_g(x,p_i) \leq r_1\}.
	\label{Eq:08VII19-B2} 
\end{equation}

Fix some $r_1 \geq s \gg i r_i$. We choose $a_i = a_{i,s} := \max_{\partial B(p_i,s)} u_i s^q$ and $b_i = b u_i(p_i)^{-\theta}$ for some large $b > 0$ (which is independent of $i$) so that, in view of \eqref{Eq:08VII19-A2}, $\varphi_i \geq u_i$ on $\partial B(p_i,s)$ and on $\partial B(p_i,ir_i)$. We then deduce from \eqref{Eq:08VII19-A2}-\eqref{Eq:08VII19-B2} and the comparison principle that
\begin{equation}
u_i \leq \varphi_i \text{ in } \{ir_i \leq d_g(x,p_i) \leq s\}.
	\label{Eq:08VII19-B3}
\end{equation}
Recalling the isolated simplicity and the gradient estimate \eqref{Eq:08VII19-A1}, we deduce from \eqref{Eq:08VII19-B3} that
\begin{align*}
s^{\frac{n-2}{2} - q}a_{i,s}
 	&= s^{\frac{n-2}{2}} \max_{\partial B(p_i,s)} u_i
		\leq Cd_g(x,p)^{\frac{n-2}{2}} \varphi_i(x) \\
	&\leq Ca_{i,s} d_g(x,p_i)^{\frac{n-2}{2} -q} + C u_i(p_i)^{-\theta} d_g(x,p_i)^{-\frac{n-2}{2} +q}
	\text{ in } \{ir_i \leq d_g(x,p_i) \leq s\}.
\end{align*}
Picking $x \in \partial B(p_i,s/C)$ for some sufficiently large $C$ and noting that $q < \frac{n-2}{2}$, we deduce that $a_{i,s} \leq C u_i(p_i)^{-\theta} s^{-(n-2-2q)}$, which gives
\[
\max_{\partial B(p_i,s)} u_i = a_{i,s} s^{-q} \leq C u_i(p_i)^{-\theta} s^{-(n-2-q)}.
\]
We have thus shown that there is some $C > 1$ so that
\[
u_i \leq C u_i(p_i)^{-\theta} s^{-(n-2-q)} \text{ in } \{Cir_i \leq d_g(x,p_i) \leq r_1\}.
\]
Estimate \eqref{Eq:07VIII19-E1} follows from the above inequality, the gradient estimate \eqref{Eq:08VII19-A1} (applied in the region $\{d_g(x,p_i) \geq r_1\}$), and \eqref{Eq:08VII19-A2} (applied in the region $\{ir_i \leq d_g(x,p_i) \leq Cr_i\}$).
\end{proof}

\begin{proof}[Proof of Lemma \ref{Lem:07VIII19-Iso+TGeo->Sim}]
Let $r_i = \varkappa^{\frac{2}{n-2}} u_i(p_i)^{-\frac{2}{n-2}}$ and
\[
\bar u_i(r) = \frac{1}{|\partial B_g(p_i, r)|_g}\int_{\partial B_g(p_i, r)} u_i(x)dS_g(x).
\]
Suppose by contradiction that $(H3)$ does not hold. Then, in view of \eqref{Eq:05VI19-A2}, there exist $\rho_i \rightarrow 0$, $\rho_i > ir_i$ such that $r^{\frac{n-2}{2}}\bar u_i(r)$ is decreasing in $(2r_i,\rho_i)$ and 
\[
\frac{d}{dr}\big|_{r = \rho_i} r^{\frac{n-2}{2}}\bar u_i(r) = 0.
\]

Define a map $\Psi_i: \RR^n \approx T_{p_i}(M,g) \rightarrow M$ by 
\[ 
\Psi_i(x) = \exp_{p_i} (\rho_i x),
\] 
and let  
\[ 
\hat v_i(x) = \rho_i^{\frac{n-2}{2}}\,u_i \circ \Psi_i(x), \qquad x  \in \RR^n. 
\] 
Then $\hat v_i$ satisfies 
\begin{equation}
f(\lambda(A_{(\hat h_i)_{\hat v_i}})) = 1  \text{ in } \{|x| < \delta_0\,\rho_i^{-1}\}, 
	\label{Eq:07VIII19-X2}
\end{equation} 
where $\hat h_i := \rho_i^{2} \Phi_i^* g$ and $\delta_0$ is the injectivity radius of $(M,g)$. Note that $\{\hat h_i\}$ converges in $C^3_{\rm loc}(\RR^n)$ to the Euclidean metric $\mathring{g}$ on $\RR^n$. Clearly,
\begin{equation}
\sup_{\{|x| < \delta_0\,\rho_i^{-1}\}} \hat v_i  
	= \hat v_i(0) 
	= \rho_i^{\frac{n-2}{2}}\,u_i(p_i) > i^{\frac{n-2}{2}} \rightarrow \infty,
	\label{Eq:07VIII19-B1}
\end{equation}
As $\{u_i\}$ is an isolated blow-up sequence, we have
\begin{equation}
\sup_{\{|x| < \delta_0\,\rho_i^{-1}\}} |x|^{\frac{n-2}{2}} \hat v_i
	= \sup_{\{d_g(x,p_i) < \delta_0\}} d_g(x,p_i)^{\frac{n-2}{2}} u_i \leq C.
	\label{Eq:07VIII19-B2}
\end{equation}
As $\{u_i\}$ has tame geometry, we also have for some $\theta \in [0,1)$ that
\begin{equation}
|\Ric_{(h_i)_{\hat v_i}}(x)|_{{(h_i)_{\hat v_i}}}
	\leq C \hat v_i(0)^{\frac{4\theta}{n-2}} |x|^{2\theta} \text{ in } \{|x| < \delta_0\,\rho_i^{-1}\}.
	\label{Eq:07VIII19-B3}
\end{equation}
Furthermore, if we let
\[
\bar v_i(r) = \frac{1}{|\partial B_{h_i}(0, r)|_{h_i}}\int_{\partial B_{h_i}(0, r)} \hat v_i(x)dS_{h_i}(x),
\]
then by contradiction hypothesis, 
\begin{equation}
r^{\frac{n-2}{2}} \bar v_i(r) \text{ is decreasing in }(2\varkappa^{\frac{2}{n-2}}\hat v_i(0)^{-\frac{2}{n-2}}, 1),
	\label{Eq:07VIII19-B4}
\end{equation}
and
\begin{equation}
\frac{d}{dr}\Big|_{r = 1}(r^{\frac{n-2}{2}} \bar v_i(r)) = 0.
	\label{Eq:07VIII19-B5}
\end{equation}

In effect, in view of \eqref{Eq:07VIII19-B1}-\eqref{Eq:07VIII19-B4}, we have rescaled $\{u_i\}$ to obtain an isolated simple blow-up sequence of solutions to \eqref{Eq:07VIII19-X2} which has tame geometry. We can then follow the proof of Lemma \ref{Lem:08VI19-AlmostSharpDecay} to show that, for any $\hat \theta \in (0,1)$,
\begin{equation}
\hat v_i(x) \leq C_{\hat\theta} \hat v_i(0)^{-\hat\theta} |x|^{-\frac{(1 + \hat \theta)(n-2)}{2}} \text{ in } \{ir_i \rho_i^{-1} \leq |x|  \leq 1\}.
	\label{Eq:07VIII19-T1}
\end{equation}

Fix some $e$ with $|e| = 1$. Define
\[
\check v_i = \frac{1}{\hat v_i(e)} \hat v_i.
\]
By \eqref{Eq:08VII19-A1}, we have
\[
|\nabla^\ell \ln \hat v_i(x)| \leq C|x|^{-\ell} \text{ in } \{0 < |x| < \delta_0\,\rho_i^{-1}\} \text{ for } \ell = 1,2.
\]
Hence, as $\check v_i(e) = 1$, $\{\check v_i\}$ converges, along a subsequence, in $C^{1,\alpha}_{loc}(\RR^n)$ to some positive function $\check v_* \in C^{1,1}_{loc}(\RR^n)$, which in view of \eqref{Eq:07VIII19-X2} and \eqref{Eq:07VIII19-T1}, satisfies
\[
\lambda(A_{\mathring{g}_{\check v_*}}) \in \partial  \Gamma \text{ in } \RR^n \setminus \{0\}.
\]
By the Liouville theorem \cite[Theorem 1.18]{Li09-CPAM} and the classification result \cite[Theorem 2.2]{LiNgBocher}, we have
\[
\check v_*(x) = \check v_*(|x|) =   (C_1 |x|^{-m} + C_2)^{\frac{n-2}{m}}   ,
\]
for some constants $m > 0$, $C_1, C_2 \geq 0$ with $C_1 + C_2 > 0$. By \eqref{Eq:07VIII19-B4}, we have that $r^{\frac{n-2}{2}} \check v_*(r)$ is decreasing in $(0,1)$ and so $C_1 > 0$. By \eqref{Eq:07VIII19-B5}, we have that $C_2 = C_1$. So 
\begin{equation}
\check v_*(x) = C_* ( |x|^{-m}+ 1)^{\frac{n-2}{m}}  \text{ for some } C_* > 0.
	\label{Eq:07VIII19T2}
\end{equation}

On the other hand, by \eqref{Eq:07VIII19-B3}, we have
\[
|\Ric_{(h_i)_{\check v_i}}(x)|_{{(h_i)_{\check v_i}}}
	\leq C\hat v_i(e)^{\frac{4}{n-2}} \hat v_i(0)^{\frac{4\theta}{n-2}} |x|^{2\theta} \text{ in } \{|x| < \delta_0\,\rho_i^{-1}\}.
\]
In view of \eqref{Eq:07VIII19-T1}, we have $\hat v_i(e) \hat v_i(0)^{\theta} \rightarrow 0$ as $i \rightarrow \infty$. This then implies (see \cite[Section 3.1, Step 6]{LiNgPoorMan}) that
\[
\Ric_{\mathring{g}_{\check v_*}} \equiv 0 \text{ in } \RR^n \setminus \{0\}.
\]
On the other hand, by \eqref{Eq:07VIII19T2}, we have
\[
R_{\mathring{g}_{\check v_*}}
	= -\frac{2}{n-2} \check v_*^{\frac{n+2}{n-2}}\Delta_{\mathring{g}} \check v_*
	\not\equiv 0,
\]
and have thus reached a contradiction.
\end{proof}

\subsubsection{Upper bound for $\tilde u_i$ and proof of Theorem \ref{Thm:deltaRHS-IsoTGeoBlUp}}\label{SSSec:ProofT1.4}

The following lemma gives the sharp upper bound for $u_i$ away from $p_\infty$. Compare Lemma \ref{Lem:08VI19-LowerSharpDecay}.

\begin{lemma}\label{Lem:08VI19-UpperSharpDecay}
Under the assumptions of Theorem \ref{Thm:deltaRHS-IsoTGeoBlUp}, for every $r \in (0,r_0/2)$, there exists $C = C(r) > 1$ (independent of $i$) such that, for all sufficiently large $i$,
\[
u_i(x) \leq C(r) u_i(p_i)^{-1} \text{ in } \{r_0/2 \geq  d_g(x,p_i) \geq r\}.
\]
\end{lemma}

\begin{proof}
Fix some $r > 0$. Suppose by contradiction that there exists $\{q_i\} \subset M$ with $r_0/2 \geq d_g(q_i,p_i) \geq r$ such that,  along a subsequence,
\begin{equation}
u_i(q_i) u_i(p_i) \rightarrow \infty.
	\label{Eq:07VIII19-SIContra}
\end{equation}

Consider the sequence
\[
\check u_i = \frac{1}{u_i(q_i)} u_i.
\]
We have $\check u_i(q_i) = 1$ and by the first and second derivative estimates \eqref{Eq:08VII19-A1}, $\{\check u_i\}$ converges, along a subsequence, in $C^{1,\alpha}_{\rm loc}( B_g(p_\infty,3r_0/4) \setminus \{p_\infty\})$ to some positive function $\check u_\infty \in C^{1,1}_{\rm loc}( B_g(p_\infty,3r_0/4) \setminus \{p_\infty\})$. By $(H2)$,
\[
|\Ric_{g_{\check u_i}}|_{g_{\check u_i}} \leq Cu_i(q_i)^{\frac{4}{n-2}} \max(1 , u_i(p_i)^{\frac{4\theta}{n-2}} d_g( \cdot, p_i)^{2\theta})  \text{ in } B_g(p_i,r_0).
\]
On the other hand, by Lemma \ref{Lem:08VI19-AlmostSharpDecay}, we have that
\[
u_i(q_i) u_i(p_i)^{\theta} \rightarrow 0 \text{ as } i \rightarrow \infty.
\]
We claim that this implies $\check u_\infty$ is smooth in $B_g(p_\infty,3r_0/4) \setminus \{p_\infty\}$ and
\begin{equation}
\Ric_{g_{\check u_\infty}} \equiv 0 \text{ in }  B_g(p_\infty,3r_0/4) \setminus \{p_\infty\}.
	\label{Eq:07VIII19-M1}
\end{equation}

Indeed, from the above, we have that
\[
-L_g \check u_i = o(1) \check u_i^{\frac{n+2}{n-2}} \text{ in } B_g(p_i,3r_0/4)
\]
where $o(1)$ denotes some function which goes to $0$ uniformly as $i \rightarrow \infty$. The convergence of $\check u_i$ to $\check u_\infty$ then implies that $\check u_\infty$ satisfies 
\[
-L_g \check u_\infty = 0 \text{ on } B_g(p_\infty,3r_0/4) \setminus \{p_\infty\} \text{ in the weak sense}.
\]
Elliptic regularity theories then imply that $\check u_\infty$ is smooth on $M\setminus \{p_\infty\}$. We can then follow \cite[Section 3.1, Step 6]{LiNgPoorMan} to obtain \eqref{Eq:07VIII19-M1}. The claim is proved.

Since $-L_g \check u_\infty = 0$ in  $B_g(p_\infty,3r_0/4) \setminus \{p_\infty\}$, we have $\check u_\infty(y) = a G_{p_\infty} + b(y)$ for some constant $a \geq 0$ and some function $b$ smooth in $B_g(p_\infty,3r_0/4)$, where $G_{p_\infty}$ is the Green's function for the conformal Laplacian with pole at $p_\infty$. By Lemma \ref{Lem:07VIII19-Iso+TGeo->Sim}, $a > 0$. Note also that, by $(H2)$,
\[
-L_{g} u_i = K_i u_i^{\frac{n+2}{n-2}},  \text{ and } |K_i| \leq C \max(1 , u_i(p_i)^{\frac{4\theta}{n-2}} d_g( \cdot, p_i)^{2\theta}) \text{ in } B_g(p_i,r_0).
\]
We now follow an argument in \cite{Li95-JDE} (see the equations (2.14)--(2.18) there) to reach a contradiction. Indeed multiplying the above equation by $u_i(q_i)^{-1}$ and integrating over a ball $B_g(p_i,r_1)$ with $0 < r_1 \ll r_0$, we get on one hand that
\begin{align*}
\limsup_{i \rightarrow 0} u_i(q_i)^{-1}\int_{B_g(p_i,r_1)} L_{g} u_i\,dx 
	&\leq \limsup_{i \rightarrow 0} u_i(q_i)^{-1}\int_{B_g(p_i,r_1)} \Delta_{g} u_i\,dx\\
	&= \int_{\partial B_g(p_\infty,r_1)} \partial_\nu (a G_{p_\infty})\,dS + O(r_1^{n-1}) < 0
\end{align*}
and on the other hand that, by Lemma \ref{Lem:08VI19-AlmostSharpDecay} and for $\frac{n+4\theta-2}{n+2} < \theta' < 1$,
\begin{align*}
&u_i(q_i)^{-1}\int_{B_g(p_i,r_1)} |K_i| u_i^{\frac{n+2}{n-2}}\,dx \\
	&\quad \leq u_i(q_i)^{-1}\int_{\{d_g(p_i,x) \leq \varkappa^{\frac{2}{n-2}} iu_i(p_i)^{-\frac{2}{n-2}}\}}   \max(1 , u_i(p_i)^{\frac{4\theta}{n-2}} d_g( \cdot, p_i)^{2\theta}) u_i^{\frac{n+2}{n-2}}\,dx \\
		&\qquad + u_i(q_i)^{-1} u_i(p_i)^{\frac{4\theta - (n+2)\theta'}{n-2}}  \int_{\{r_1 \geq d_g(p_i,x) \geq  \varkappa^{\frac{2}{n-2}} i u_i(p_i)^{-\frac{2}{n-2}}\}} d_g(x,p_i)^{2\theta -\frac{(1 + \theta')(n+2)}{2}}\,dx \\
	&\quad \leq u_i(q_i)^{-1} u_i(p_i)^{-1} O(1)
	\rightarrow 0
\end{align*}
which amounts to a contradiction.
\end{proof}

\begin{proof}[Proof of Theorem \ref{Thm:deltaRHS-IsoTGeoBlUp}]
By Lemma \ref{Lem:08VI19-UpperSharpDecay}, $\{\tilde u_i\}$ is bounded in $C^0_{loc}( B_g(p_\infty,r_0/2) \setminus \{p_\infty\})$. By estimate \eqref{Eq:08VII19-A1}, $\{\tilde u_i\}$ converges along a subsequence in $C^{1,\alpha}_{\rm loc}( B_g(p_\infty,r_0/2) \setminus \{p_\infty\})$ to some positive function $\tilde u_\infty \in C^{1,1}_{\rm loc}( B_g(p_\infty,r_0/2) \setminus \{p_\infty\})$. Moreover, the same argument giving \eqref{Eq:07VIII19-M1} shows that $\tilde u_\infty \in C^{\infty}(B_g(p_\infty,r_0/2) \setminus \{p_\infty\})$ and
\begin{equation}
A_{g_{\tilde u_\infty}} \equiv 0  \text{ in } B_g(p_\infty,r_0/2) \setminus \{p_\infty\}.
	\label{Eq:07VIII19-Y1}
\end{equation}

We claim that 
\[
c := \lim_{d_g(x,p_\infty) \rightarrow 0} \tilde u_\infty(x) d_g(x,p_\infty)^{n-2} \text{ exists and is positive}.
\]
First, by Lemma \ref{Lem:08VI19-LowerSharpDecay}, 
\[
\underline{c} := \liminf_{d_g(x,p_\infty) \rightarrow 0} \tilde u_\infty(x) d_g(x,p_\infty)^{n-2} \text{ is finite and positive.}
\]
The claim is then proved by following Step 4 in the proof of \cite[Theorem 1.3]{LiNgPoorMan}, which we briefly outline here for readers' convenience. By \eqref{Eq:08VII19-A1}, this implies that
\[
\overline{c} := \limsup_{d_g(x,p_\infty) \rightarrow 0} \tilde u_\infty(x) d_g(x,p_\infty)^{n-2} \text{ is finite and non-negative.}
\]
Now if $\underline{c} < \overline{c}$, then by performing a blow-up argument at $p_\infty$, we would obtain a function $\tilde v \in C^{1,1}_{\rm loc}(\RR^n \setminus \{0\})$ satisfying $\lambda(A_{\mathring{g}_{\check v}}) \in \partial \Gamma_k$ in $\RR^n \setminus \{0\}$ and
\[
\min_{|x| =1} \tilde v(x) < \sup_{|x| = 1} \tilde v(x)
\]
which would contradict the symmetry result \cite[Theorem 1.18]{Li09-CPAM}. We conclude that $\underline{c} = \overline{c}$ and so $c$ exists as desired; see \cite{LiNgPoorMan} for details.

By Corollary \ref{Cor:sigmakL1Bnd}, the restriction of $\tilde u_i^{\frac{n+2k}{n-2k}} \sigma_k(A_{g_{\tilde u_i}})$ to $B(p_\infty,r_0/2)$ weakly* converges to $m_{n,k}\delta_{p_\infty}$ with (see also Proposition \ref{Prop:CompCptness1})
\[
m_{n,k} = U_1(0)^{\frac{n-2k}{n-2}} \int_{\RR^n} U_1^{\frac{n+2k}{n-2}}\,dx. 
\]
The proof is complete.
\end{proof}

\subsection{A divergence identity and its consequences}\label{SSec:DivId}

In this subsection, we present a divergence identity for the Newton tensors associated with the $(1,1)$-Schouten tensor. 

For a symmetric $(1,1)$-tensor $A$, the symmetric functions $\sigma_0(A), \ldots, \sigma_n(A)$ are defined by
\[
\det (\lambda I - A) = \sum_{k = 0}^n (-1)^k \sigma_k(A)\,\lambda^{n-k}.
\]
It is clear that $\sigma_k(A) = \sigma_k(\lambda(A))$. The Newton tensors $\NewtonT{k}(A)$, $k = 0, \ldots, n-1$, of $A$ are defined by
\[
(\NewtonT{k}(A))^i{}_j = \frac{\partial \sigma_{k+1}}{\partial A^j{}_i}(A).
\]
It is well known that
\[
\NewtonT{k}(A) = \sum_{l = 0}^k (-1)^{k-l}\,\sigma_l(A)\, A^{k-l},
\]
and, for $0 \leq k \leq n -1$,
\begin{align}
\tr\,\NewtonT{k}(A) 
	&= (n-k)\sigma_k(A) ,\label{Eq:TTrace}\\
\NewtonT{k+1}(A) 
	&= - A \NewtonT{k}(A) +  \sigma_{k+1}(A)\,I ,
	\label{Eq:TRecursion}\\
\tr(A \NewtonT{k}(A))
	&= (k+1)\sigma_{k+1}(A).
	\label{Eq:trTA}
\end{align}

In the sequel, for a given metric $g$, we use $\NewtonT{k}(A_g)$ to denote the Newton tensors of the $(1,1)$-Schouten tensor $A_g$.

When $g$ is locally conformally flat, it is well known that $\NewtonT{k}(A_{g_u})$ has a divergence structure, see \cite{Viac00-Duke}. The following lemma gives a generalization of that statement.

\begin{lemma}\label{Lem:DivDF}
Let $U$ be an $n$-dimensional manifold with or without boundary, $g$ be a smooth Riemannian metric on $U$, and let $0 \leq k \leq n-1$. For any smooth positive function $u$ on $U$, we have
\begin{align}
\nabla_j \NewtonT{k}{}^j{}_r(A_{g_u}) 
	&= -\frac{2}{n-2} \frac{\nabla_j u}{u}\Big[n\,\NewtonT{k}{}^j{}_r(A_{g_u}) - (n-k)\sigma_k(A_{g_u})\delta^j{}_{r}\Big]\nonumber\\
		&\quad
			+ \frac{1}{n-2} u^{-\frac{n+2}{n-2}}  \sum_{q=1}^{k-1} (-1)^{k-q} \NewtonT{q}{}^j{}_l(A_{g_u}) [-2 W_{jt}{}^l{}_s \nabla^s u + u C^l{}_{tj}] (A_{g_u}^{k-1-q})^t{}_r,
	\label{Eq:12III19-DivId}
\end{align}
where $\nabla$ is the covariant derivative of $g$, and $W$ and $C$ are the Weyl and Cotton tensors of $g$, and, for $k = 0$ or $k = 1$, the summation on the right hand side is trivial.
\end{lemma}

\begin{remark}
If we let $\widetilde \nabla$ denote the covariant derivative of $g_u$, then \eqref{Eq:12III19-DivId} is equivalent to
\begin{align}
\widetilde \nabla_j \NewtonT{k}{}^j{}_r(A_{g_u}) 
	&= - \frac{1}{n-2} u^{-\frac{n+2}{n-2}}  \sum_{q=1}^{k-1} (-1)^{k-q} \NewtonT{q}{}^j{}_l(A_{g_u}) [-2 W_{jt}{}^l{}_s \nabla^s u + u C^l{}_{tj}] (A_{g_u}^{k-1-q})^t{}_r.
	\label{Eq:12III19-DivIdX}
\end{align}
(In particular, if $g$ is locally conformally flat or $k = 0$ or $k = 1$, $\NewtonT{k}(A_{g_u})$ is divergence-free with respect to $g_u$.) Similarly, identity \eqref{Eq:Acom-1} below is equivalent to
\begin{align}
\widetilde\nabla_i (A_{g_u})^l{}_j - \widetilde\nabla_j (A_{g_u})^l{}_i
	=  \frac{1}{2}u^{-\frac{n+2}{n-2}}\,W_{ij}{}^l{}_{s}\,\nabla^s u + \frac{1}{n-2} u^{-\frac{4}{n-2}}\,C^l{}_{ji}.
			\label{Eq:Acom-1X}
\end{align}
In view of the identity $\sigma_k(A) = \frac{1}{k} \tr\,(\NewtonT{k-1}(A) A)$, the identities \eqref{Eq:12III19-DivIdX} and \eqref{Eq:Acom-1X} give a div-curl structure for the $\sigma_k$ operator.
\end{remark}

As an application of Lemma \ref{Lem:DivDF}, we establish the following compensated compactness result for the $\sigma_k$ equation.

\begin{proposition}\label{Prop:CompCptness1}
Let $U$ be a compact $n$-dimensional manifold with or without boundary, $g$ be a smooth Riemannian metric on $U$, and let $1 \leq k \leq n$. Suppose $\{u_j\}$ is a sequence of smooth positive functions on $U$ which converges in $C^{1,\alpha}(U,g)$ for some $0 < \alpha < 1$ and weakly in $W^{2,k}(U,g)$ to some positive function $u \in C^{1,\alpha}(U,g) \cap W^{2,k}(U,g)$. Then, for all $\varphi \in C^0(U)$ satisfying $\varphi = 0$ on $\partial U$,
\[
\lim_{j \rightarrow \infty} \int_{U}  \sigma_k (\lambda (A_{g_{u_j}}))\varphi \,dv_{g} =  \int_{U} \sigma_k (\lambda (A_{g_{u}}))\varphi \,dv_{g}.
 \]
\end{proposition}

\begin{corollary}\label{Cor:CC1}
Under the assumptions of Proposition \ref{Prop:CompCptness1}, one has for all $\gamma \in \RR$ and $\varphi \in C^0(U)$ satisfying $\varphi = 0$ on $\partial U$ that
\[
\lim_{j \rightarrow \infty} \int_{U} u_j^{\gamma} \sigma_k (\lambda (A_{g_{u_j}}))\varphi \,dv_{g} =  \int_{U} u^{\gamma} \sigma_k (\lambda (A_{g_{u}}))\varphi \,dv_{g}.
 \]
\end{corollary}

\begin{proof}[Proof of Lemma \ref{Lem:DivDF}]
It is more convenient to work with $w = u^{-\frac{2}{n-2}}$ so that
the $(1,1)$-Schouten tensor of $g_u = w^{-2} g$ is given by
\[
(A_w)^i{}_j  = w\nabla^i \nabla_j  w - \frac{1}{2} |\nabla w|_{g}^2\,\delta^i{}_j + w^2 A^i{}_j,
\]
where $A = A_g$ is the $(1,1)$-Schouten tensor of $g$.
 
In the proof, indices are lowered and raised using $g$.

 Fix a point $p$ and let $\{x^1, \ldots, x^n\}$ be a geodesic normal coordinate system at $p$. In particular, $\Gamma_{ij}^l(p) = 0$. The following computation is done at $p$.

First, we have
\begin{align}
\nabla_i (A_w)^l{}_j - \nabla_j (A_w)^l{}_i
	&= \nabla_i w \nabla^l \nabla_{j} w - \nabla_j w \nabla^l\nabla_{i} w - \nabla_{is} w\,\nabla^s w\,\delta^l{}_{j} + \nabla_{js} w\,\nabla^s w\,\delta^l{}_{i}\nonumber\\
		&\quad +w\,\textrm{Riem}_{ij}{}^l{}_{s}\,\nabla^s w
			 + 2w\,(\nabla_i w\,A^l{}_j - \nabla_j w\,A^l{}_i) + \frac{1}{n-2} w^2\,C^l{}_{ji}\nonumber\\
	&= \frac{1}{w}\Big[\nabla_i  w\,(A_w)^l{}_j  - \nabla_j  w\,(A_w)^l{}_i \nonumber\\
			&\qquad - \nabla_s w (A_w)^s{}_{i}\,\delta^l{}_j + \nabla_s w (A_w)^s{}_{j}\,\delta^l{}_i\Big]\nonumber\\
		&\quad + w\,W_{ij}{}^l{}_{s}\,\nabla^s w + \frac{1}{n-2} w^2\,C^l{}_{ji},
			\label{Eq:Acom-1}
\end{align}
where $\textrm{Riem}$ is the Riemann curvature tensor of $g$ and where we have used the Ricci decomposition $\textrm{Riem}_{ij}{}^l{}_s
	=  - A^l{}_j g_{is} + A^l{}_i g_{js}  - A_{is} \delta^l{}_j  + A_{js} \delta^l{}_i+ W_{ij}{}^l{}_s$.

Using \eqref{Eq:TTrace}, \eqref{Eq:TRecursion} and \eqref{Eq:Acom-1}, we compute
\begin{eqnarray}
\nabla_j \NewtonT{k+1}{}^j{}_r (A_w)
	&\stackrel{\eqref{Eq:TRecursion} }{=}&
	 - \nabla_j \NewtonT{k}{}^j{}_l(A_w) (A_w)^l{}_r -  \NewtonT{k}{}^j{}_l(A_w) [\nabla_j (A_w)^l{}_r - \nabla_r (A_w)^l{}_j]\nonumber\\
	&\stackrel{\eqref{Eq:Acom-1}}{=}& - \nabla_j \NewtonT{k}{}^j{}_l(A_w) (A_w)^l{}_r \nonumber\\
		&&\quad -  \NewtonT{k}{}^j{}_l(A_w) \Big\{\frac{1}{w}\big[\nabla_j  w\,(A_w)^l{}_r  - \nabla_r  w\,(A_w)^l{}_j \nonumber\\
		&&\quad - \nabla_s w (A_w)^s{}_j\,\delta^l{}_r + \nabla_s w (A_w)^s{}_r\,\delta^l{}_j\big]
		\nonumber \\
		&&\quad
			+ w\,W_{jr}{}^l{}_{s}\,\nabla^s w + \frac{1}{n-2}w^2\,C^l{}_{rj}
		\Big\}\nonumber\\
	&\stackrel{\eqref{Eq:TTrace},\eqref{Eq:TRecursion} }{=}&
	- \nabla_j \NewtonT{k}{}^j{}_l(A_w) (A_w)^l{}_r 
		+ \frac{k+1}{w} \nabla_r w \sigma_{k+1}(A_w)\nonumber\\
		&&\quad
		 -\frac{n-k}{w} \nabla_s w \sigma_k(A_w) (A_w)^s{}_r \nonumber\\
		&&\quad 
			-  w\NewtonT{k}{}^j{}_l(A_w)\big[ W_{jr}{}^l{}_{s}\,\nabla^s w + \frac{1}{n-2}w\,C^l{}_{rj}
		\big].
		\label{Eq:12III19-V2}
\end{eqnarray}

Identity \eqref{Eq:12III19-DivId} then follows from an induction on $k$ using \eqref{Eq:12III19-V2}.
\end{proof}

\begin{proof}[Proof of Proposition \ref{Prop:CompCptness1}] The result is clear for $k = 1$. Suppose that $2 \leq k \leq n$.  

Using a partition of unity if necessary, we may assume for simplicity that $U$ is contained in a single chart. 

Let $A_j$ and $A$ denote the $(1,1)$-Schouten tensor of $g_{u_j}$ and $g_u$, and $\NewtonT{\ell}_j$ and $\NewtonT{\ell}$ denote the $\ell$-th Newton tensor of $A_j$ or $A$, respectively.

By the hypotheses, $A_j$ converges weakly in $L^k(U,g)$ to $A$. Also, for $1 \leq \ell \leq k - 1$, $\{\NewtonT{\ell}_j\}$ is bounded in $L^{k/\ell}(U,g)$ and so converges weakly in $L^{k/\ell}(U,g)$ to some $\NewtonT{\ell}_\infty$. 

We first show that $\NewtonT{\ell}_\infty = \NewtonT{\ell}$ for $1 \leq \ell \leq k - 1$ by an induction on $\ell$. For $\ell = 1$, the assertion holds due to the weak convergence of $A_j$ to $A$. Assume that the assertion holds for some $\ell \leq k - 2$. Recall that, by Lemma \ref{Lem:DivDF}, the divergence of each column of $\NewtonT{\ell}_j$ is bounded in $L^{k/\ell}(U,g)$, and by \eqref{Eq:Acom-1}, the curl of each row of $A_j$ is bounded in $L^k(U,g)$. An application of the div-curl lemma (\cite{Murat78}) then implies that $\{A_j \NewtonT{\ell}_j\}$ converges to $A\,\NewtonT{\ell} $ in the sense of distribution. In view of \eqref{Eq:TRecursion}-\eqref{Eq:trTA}, this implies that $\{\NewtonT{\ell+1}_j\}$ converges to $\NewtonT{\ell + 1}$ in the sense of distribution, from which we conclude that $\NewtonT{\ell+1}_\infty = \NewtonT{\ell + 1}$.

The argument above in fact also shows that $\{A_j \NewtonT{k - 1}_j\}$ converges to $A\,\NewtonT{k - 1}$ in the sense of distribution. By \eqref{Eq:trTA}, this implies that $\{\sigma_k(A_j)\}$ converges to $\sigma_k(A)$ in the sense of distribution. Recalling that $\{\sigma_k(A_j)\}$ is bounded in $L^1(U,g)$, we are done.
\end{proof}


\subsection*{Acknowledgement.} Both authors would like to thank Xiaochun Rong and Jiaping Wang for stimulating discussions. Y.Y. Li is partially supported by NSF grant DMS-1501004.

\appendix

\section{Smooth concave defining functions of cones}\label{Sec:SDefFc}

In this appendix, we construct for every given $\Gamma$ satisfying \eqref{G1}-\eqref{G2} a function $f$ satisfying  \eqref{f1}-\eqref{f3}, which was used in the proof of Theorem \ref{thm:MainThm}(i).

\begin{proposition}\label{prop:fConstr}
Let $\Gamma$ satisfy \eqref{G1}-\eqref{G2}. Then there exists a concave function $f \in C^\infty(\Gamma) \cap C(\bar \Gamma)$ satisfying \eqref{f1}-\eqref{f3}. If it holds in addition that $(1,0,\ldots, 0) \in \Gamma$, then there exists $\nu \in ( 0,1)$ such that
\begin{equation}
\frac{\partial f}{\partial \lambda_i}(\lambda) \geq \nu \sum_j \frac{\partial f}{\partial \lambda_j}(\lambda) \text{ for all } i = 1, \ldots, n \text{ and } \lambda \in \Gamma.
	\label{Eq:22I19-UE}
\end{equation}
\end{proposition}

We note that condition \eqref{Eq:22I19-UE} is related to the strict ellipticity of equation \eqref{Eq:04VI19-NE}.

\begin{proof} If $\Gamma = \Gamma_1$, the result is obvious. We assume that $\Gamma \neq \Gamma_1$. Then the set $\Omega_\Gamma= \Gamma \cap \{\lambda: [\lambda] := \lambda_1 + \ldots + \lambda_n = 1\}$ is bounded and convex. It is well known that $\Omega_\Gamma$ admits a concave defining functions $h$ such that $h > 0$ in $\Omega_\Gamma$ and $h = 0$ on $\partial \Omega_\Gamma$ (see e.g. \cite[Section 2.1]{HormanderConvexity}). Furthermore, $h$ can be chosen in $C^\infty(\Omega_\Gamma) \cap C(\bar \Omega_\Gamma)$ (see e.g. \cite[Theorem 7]{ChengYau77-CPAM}). (In fact one can have $h \in C^{\beta}(\bar\Omega_\Gamma)$ with $\beta = \frac{2}{n}$ if $n \geq 3$ and $0 < \beta < 1$ if $n = 2$, but this is not needed in the present argument; see \cite[Lemma 1]{Caffarelli90-AoM}.)

 By considering 
\[
\tilde h(\lambda) = \sum_{x \text{ is a permutation of } \lambda} h(x),
\]
instead of $h$, we can assume without loss of generality that $h$ is symmetric.

 Let $\nabla_T$ denote the gradient on $\Omega_\Gamma$. Observe that for $x$ $\in$ $\Omega_\Gamma$ and $p_0$ $\in$ $\bar\Omega_\Gamma$, the concavity of $h$ implies that
\begin{equation}
h(x)  -\nabla_T h(x) \cdot (x - p_0) \geq h(p_0) \geq 0.
	\label{Eq:23II19-R1}
\end{equation}
Let
\[
\alpha = \left\{\begin{array}{cl}
\text{ any number in $(0,1)$} & \text{ if } (1, 0, \ldots, 0) \in \partial\Gamma,\\
1 & \text{ if } (1,0, \ldots, 0) \in \Gamma,
\end{array}\right.
\]
and $g$ $=$ $h^\alpha$. By \eqref{Eq:23II19-R1}, we have
\begin{multline}
g(x) - \nabla_T g(x) \cdot (x - p_0) = h(x)^{\alpha-1}\big[h(x) - \alpha\,\nabla_T h(x) \cdot (x - p_0)\big]\\
 \geq h(x)^{\alpha-1}[(1-\alpha)h(x) + \alpha h(p_0)] \text{ for any } x \in \Omega_\Gamma \text{ and } p_0 \in \bar\Omega_\Gamma.
\label{EllConstr}
\end{multline}
Note that the right hand side of \eqref{EllConstr} is non-negative and is zero if and only if $\alpha = 1$ (i.e. $(1,0, \ldots, 0) \in \Gamma$) and $p_0 \in \partial \Omega_\Gamma$.

Define $f$ by
\[
f(\lambda) = (\lambda_1 + \ldots + \lambda_n)\,g\big(\frac{\lambda}{\lambda_1 + \ldots + \lambda_n}\big).
\]
We now show that $\partial_i f$ $>$ $0$ and $f$ is concave in $\Gamma$. 

Let
\[
[\lambda] = \lambda_1 + \ldots + \lambda_n \text{ and }\lambda' = \frac{\lambda}{[\lambda]}.
\]
We compute
\begin{multline}
\partial_i f(\lambda)
	= g(\lambda') + [\lambda]\,\partial_j g(\lambda')\,\frac{\delta_{ij}[\lambda] - \lambda_j}{[\lambda]^2}\\
		= g(\lambda') + \partial_i g(\lambda') - \partial_j g(\lambda')\lambda'_j
			= g(\lambda') - \nabla_T g(\lambda') \cdot(\lambda' - p^i),
			\label{Eq:23II19-R3}
\end{multline}
where $p^i_j$ $=$ $\delta^i_j$. Since $\Gamma$ $\supset$ $\Gamma_n$, it follows that $p^i$ $\in$ $\bar\Omega_\Gamma$. Hence, by \eqref{EllConstr},
\begin{equation}
\partial_i f(\lambda) \geq h(\lambda')^{\alpha - 1}[(1-\alpha)h(\lambda') + \alpha h(p^i)] \text{ in }\Gamma.
	\label{Eq:23II19-R2}
\end{equation}
If $(1,0, \ldots, 0) \in \partial\Gamma$, then $\alpha \in (0,1)$ and so the right hand side of \eqref{Eq:23II19-R2} is larger or equal to $(1-\alpha)h(\lambda')^\alpha > 0$. If $(1,0, \ldots, 0) \in \Gamma$, then $p^i \in \Omega_\Gamma$ and so the right hand side of \eqref{Eq:23II19-R2} is larger or equal to $h(\lambda')^{\alpha - 1} h(p^i) > 0$. In either case, we have
\[
\partial_i f(\lambda)  > 0 \text{ in }\Gamma.
\]

To prove the concavity of $f$, we calculate its Hessian. We have
\begin{align*}
[\lambda]\partial_{ij} f(\lambda)
	&= \partial_k g(\lambda')\frac{\delta_{kj}[\lambda] - \lambda_k}{[\lambda]} + \partial_{ki} g(\lambda')\frac{\delta_{kj}[\lambda] - \lambda_k}{[\lambda]}\\
		&\qquad\qquad - \partial_l g(\lambda')\frac{\delta_{lj}[\lambda] - \lambda_l}{[\lambda]} - \partial_{kl} g(\lambda')\,\lambda'_l\,\frac{\delta_{kj}[\lambda] - \lambda_k}{[\lambda]}\\
	&= \partial_{ij} g(\lambda') - \partial_{ki} g(\lambda')\lambda'_k - \partial_{lj} g(\lambda')\lambda'_l + \partial_{kl} g(\lambda')\lambda'_k\,\lambda'_l.
\end{align*}
Hence, for any $p$ $\in$ $\RR^n$, we have
\begin{align*}
[\lambda]\partial_{ij} f(\lambda)\,p_i\,p_j
	&= \partial_{ij} g(\lambda')\,p_i\,p_j - \partial_{ki} g(\lambda')\lambda'_k\,p_i\,p_j\\
		&\qquad\qquad - \partial_{lj} g(\lambda')\lambda'_l\,p_i\,p_j + \partial_{kl} g(\lambda')\lambda'_k\,\lambda'_l\,p_i\,p_j\\
	&= \partial_{ij} g(\lambda')\,p_i\,p_j - 2\partial_{ki} g(\lambda')\lambda'_k\,p_i\,[p] + \partial_{kl} g(\lambda')\lambda'_k\,\lambda'_l\,[p]^2\\
	&= \partial_{ij} g(\lambda')\,(p_i - \lambda'_j[p])(p_j - \lambda'_j[p]) \leq 0,
\end{align*}
where we have used $\nabla_T^2 g$ $\leq$ $0$ in $\Omega_\Gamma$. As $\Gamma$ $\subset$ $\Gamma_1$, $[\lambda]$ $>$ $0$ in $\Gamma$. Therefore, $\nabla^2 f$ $\leq$ $0$ in $\Gamma$, i.e. $f$ is concave in $\Gamma$.

Finally, assume that $(1,0, \ldots, 0)$ is in $\Gamma$, we show that \eqref{Eq:22I19-UE} holds. For any $x \in \Omega_\Gamma$, define $L_x: \Omega_\Gamma \rightarrow \RR$ by
\[
L_{x}(p) = g(x) - \nabla_T g(x) \cdot (x - p) = h(x) - \nabla_T h(x) \cdot (x - p), \qquad p \in \Omega_\Gamma.
\]
Note that $L_x$ is a linear function, and hence is harmonic with respect to the metric induced on $\Omega_\Gamma$ by the Euclidean metric on $\RR^n$. Furthermore, by \eqref{Eq:23II19-R1}, $L_x$ is positive in $\Omega_\Gamma$. Since all $p^1, \ldots, p^n \in \Omega_\Gamma$, it follows from the Harnack inequality that there is some constant $C$ depending only on $\Omega_\Gamma$ such that
\[
L_x(p^i) \leq C L_x(p^j) \text{ for all } x \in \Omega_\Gamma, 1 \leq i \leq j \leq n.
\]
Recalling \eqref{Eq:23II19-R3}, we obtain that
\[
0 < \partial_i f(\lambda) \leq C \partial_j f(\lambda) \text{ for all } \lambda \in \Gamma, 1 \leq i \leq j \leq n,
\]
which implies \eqref{Eq:22I19-UE}.
\end{proof}

\begin{proposition}\label{prop:Nof6.19}
Let $\Gamma$ satisfy \eqref{G1}-\eqref{G2}. If $(1,0,\ldots, 0) \in \partial \Gamma$, then there is no function $f \in C^\infty(\Gamma) \cap C(\bar \Gamma)$ satisfying simultaneously \eqref{f1}-\eqref{f4} and \eqref{Eq:22I19-UE}.
\end{proposition}

\begin{proof} Suppose by contradiction that there is some $f \in C^\infty(\Gamma) \cap C(\bar \Gamma)$ satisfying simultaneously \eqref{f1}-\eqref{f4} and \eqref{Eq:22I19-UE}. By \eqref{f3} and \eqref{Eq:22I19-UE}, it follows that there is some constant $C > 0$ such that
\begin{equation}
0 < \partial_i f(\lambda) \leq C \partial_j f(\lambda) \qquad \text{ for all } \lambda \in \Gamma, 1 \leq i, j \leq n.
	\label{Eq:24II19-X1}
\end{equation}

Let $\Omega_\Gamma= \Gamma \cap \{\lambda: [\lambda] := \lambda_1 + \ldots + \lambda_n = 1\}$ and $\nabla_T$ denote the gradient on $\Omega_\Gamma$. Then $\Omega_\Gamma$ is bounded and convex and $h := f|_{\Omega_\Gamma}$ is a positive concave defining function for $\Omega_\Gamma$.

We write $[\lambda] = \lambda_1 + \ldots + \lambda_n$  and $\lambda' = \frac{\lambda}{[\lambda]}$. Then with $p^i_j$ $=$ $\delta^i_j$ (see \eqref{Eq:23II19-R3}), 
\begin{equation}
\partial_i f(\lambda)
			= h(\lambda') - \nabla_T h(\lambda') \cdot(\lambda' - p^i).
			\label{Eq:24II19-X2}
\end{equation}

For any $x \in \Omega_\Gamma$, define $L_x: \Omega_\Gamma \rightarrow \RR$ by
\[
L_{x}(p) = h(x) - \nabla_T h(x) \cdot (x - p) = h(x) - \nabla_T h(x) \cdot (x - p), \qquad p \in \Omega_\Gamma.
\]
By \eqref{Eq:24II19-X1}-\eqref{Eq:24II19-X2}, we have that
\[
0 < L_x(p^i) \leq C L_x(p^j) \text{ for all } x \in \Omega_\Gamma, 1 \leq i, j \leq n.
\]
In particular, since $L_x$ is a linear function, we have that 
\[
0 < L_x(p) \leq C L_x(q) \text{ for all } x \in \Omega_\Gamma, p, q \in \Omega_n,
\]
where $\Omega_n \subset \Omega_\Gamma$ is the interior of the convex hull of the points $p^1, \ldots, p^n$. In particular, we have
\begin{equation}
0 < L_x(\frac{1}{n}(1, \ldots, 1)) =  L_x(\frac{1}{n}(p^1 + \ldots + p^n)) \leq C L_x(x) = cg(x) \text{ for all } x \in \Omega_n.
	\label{Eq:24II19-X3}
\end{equation}
On the other hand, by the concavity of $h$ on $\Omega_\Gamma$ and the definition of $L_x$, we have
\[
L_x(p) \geq h(p) \text{ for all } x, p \in \Omega_\Gamma.
\]
It follows that $L_x(\frac{1}{n}(1, \ldots, 1)) \geq h(\frac{1}{n}(1, \ldots, 1)) > 0$. Returning to \eqref{Eq:24II19-X3}, we obtain
\[
0 < h(\frac{1}{n}(1, \ldots, 1)) \leq ch(x) \text{ for all } x \in \Omega_n.
\]
Sending $x \rightarrow p^1$ for example, this implies that 
\[
0 < h(\frac{1}{n}(1, \ldots, 1)) \leq 0,
\]
which is absurd. The proposition is proved.
\end{proof}


\section{Convexity of sets of symmetric matrices and sets of eigenvalues}

We give a presumably well-known statement on eigenvalues of sums of matrices which is used in the body of the paper.
\begin{lemma}\label{Lem:CoConeMatrix}
Let $G \subset \RR^n$ be a symmetric subset of $\RR^n$ and $U \subset Sym^{n \times n}$ be the set of real symmetric $n \times n$ matrices whose eigenvalues belong to $G$. Then $G$ is convex if and only if $U$ is convex.
\end{lemma}

\begin{proof} It is clear that $G$ is convex if $U$ is convex. To prove the converse, it suffices to show that, for any symmetric matrices $A$ and $B$ with eigenvalues $u$ and $v$ respectively, the eigenvalues $w$ of $\frac{1}{2}(A + B)$ belongs to the convex hull of the set $X$ consisting of the permutations of $u$ and $v$.

Note that there exist orthogonal matrices $P$ and $Q$ such that
\begin{equation}
w_i = \frac{1}{2} \sum_{j=1}^n (P_{ij}^2\,u_j + Q_{ij}^2\,v_j), \qquad i = 1, 2, \ldots, n.
	\label{Eq:Appwi}
\end{equation}

Consider the matrix $S$ defined by $S_{ij} = P_{ij}^2$. As $P$ is orthogonal, $S$ is doubly stochastic (i.e. the entries of $S$ are non-negative and each of its rows and columns sums to one), and hence by the Birkhoff-von Neumann theorem, $S$ is a linear combination of permutation matrices. It follows that the vector $Su$ belongs to the convex hull of the permutations of $u$.

Noting that $(Su)_i = \sum_j P_{ij}^2\,u_j$, we deduce from the foregoing paragraph and \eqref{Eq:Appwi} that $w$ belongs to the convex hull of $X$, as desired.
\end{proof}



\newcommand{\noopsort}[1]{}

\end{document}